%Paper: alg-geom/9502005
%From: Igor Dolgachev <Igor.Dolgachev@math.lsa.umich.edu>
%Date: Tue, 7 Feb 1995 12:18:03 -0500
%Date (revised): Sun, 7 Jan 96 09:18:02 EST

%12/31/95
%PLAIN TEX

\magnification = 1200
\def\pol{$M$-polarized }

\def\bfr{{\bf R}}
\def\bfc{{\bf C}}
\def\bfz{{\bf Z}}
\def\bfp{{\bf P}}
\def\cf{{\cal F}}
\def\pic{{\rm Pic}}
%\font\smalltype=cmr10 at 10truept

\vglue .8 in
\centerline{{\bf MIRROR SYMMETRY FOR LATTICE POLARIZED K3 SURFACES}}

\vglue .2 in

\centerline{{\bf Igor V. Dolgachev}\footnote{$\ ^*$}{
Research supported in part by a NSF grant.}}
\bigskip
\centerline {{\it Department of Mathematics, University of Michigan, Ann Arbor,
MI 48109}}

\vglue .5 in

{\bf Introduction.}
There has been a recent explosion in the number of mathematical publications
due to the discovery
of a certain duality between some families of
Calabi-Yau threefolds made by a group of theoretical physicists (see {\bf
[11,26]}
for references). Roughly speaking this duality, called
 mirror symmetry, pairs two families ${\cal F}$ and ${\cal F}^*$ of Calabi-Yau
threefolds in such a way that the following properties are satisfied:

{\parindent=1.5 cm
\item{MS1} The choice of the mirror family ${\cal F}^*$ involves the choice of
a boundary
point $\infty$ of a compactification $\bar {\cal F}$ of the moduli space for
$\cf$ at which
the monodromy is ``maximally unipotent''.
\item{MS2} For each $V\in \cf$ and $V'\in \cf^*$ the Hodge numbers satisfy
$$h^{1,1}(V) = h^{2,1}(V'),\ \ h^{2,1}(V) = h^{1,1}(V').$$
\item{MS3} For some open subset $U$ of $\infty$, for
any $V\in U\cap \cf$, the Laurent expansion of the canonical symmetric
trilinear form
$S^3(H^1(V,\Theta_V))\to H^{0,3}(V)^{\otimes 2}$ (the Griffiths-Yukawa cubic)
at $\infty$ can be
identified, after some special choice
of local parameters and a basis of $H^{0,3}(V)$,
with the quantum intersection form on $H^2(V'), V'\in \cf^*.$
\item{MS4} The period map induces a holomorphic multivalued mapping from the
subset $U\cap \cf$
to the tube domain
$H^2(V',\bfr)+i{\cal K}_{V'}$ where ${\cal K}_{V'}$ is the K\"ahler cone of
$V'\in \cf^*$ (the mirror mapping).
\par}

Although known to some experts but never stated explicitly, it is a
fact that
 mirror symmetry is a very beautiful and non-trivial (in many
respects still hypothetical) generalization to the next dimension of the
duality for
K3 surfaces discovered almost 20 years
ago by H. Pinkham {\bf [34]} and independently by the author and V. Nikulin
{\bf
[8,9,31]}. This duality was used
to explain Arnold's Strange Duality  for exceptional unimodal critical points
{\bf [1]}.
There are repeated hints on the relationship between the latter duality and the
mirror symmetry both in
physics literature ({\bf [2,14,23]}) and in mathematics literature ({\bf
[6,19,35,41]}).
Some of the results of this paper were independently obtained in
{\bf [3,19,21,27,35]} and some must be known
to V. Batyrev and V. Nikulin. The paper {\bf [40]} of
Todorov is probably most relevant.
Nevertheless I believe that it is worthwhile to give a detailed account of how
the ideas of
Arnold's Strange Duality allows one to state (and prove) precise
analogs of properties MS1-MS4 for K3 surfaces.

Note that property MS2 says that the local moduli number of $V\in\cf$ is equal
to the second
Betti number of $V'\in\cf^*$. In the case of K3 surfaces, the first number is
always equal to 20,
and the second number is equal to 22. The key observation is that in
three-dimensional case the second Betti number
is  equal to the rank of the Picard group of algebraic cycles. This suggests to
create different
moduli families of K3 surfaces with condition on the Picard group. The simplest
realization of this idea is
based on
the notion of a polarized K3 surface.
This is pair $(X,h)$ where $X$ is a K3 surface and $h\in {\rm Pic}(X)$ is an
ample (or pseudo-ample) divisor class. A generalization of this notion,
due to V. Nikulin {\bf [30]}, is the notion of a lattice polarized K3 surface.
We fix a lattice $M$ (a free abelian group equipped with an integral quadratic
form) and consider a pair $(X,j)$ where
$X$ is a K3 surface, and $j:M\to \pic(X)$ is a primitive embedding of lattices
such that $j(M)$
contains a pseudo-ample divisor class. One can construct a coarse moduli space
${\bf K}_M$ of $M$-polarized K3 surfaces.
An obvious condition of its non-emptiness is that $M$ is isomorpic to a
sublattice of an even
unimodular lattice $L$ of signature $(3,19)$ isomorphic to the second
cohomology group of a K3 surface equipped with the cup-product.

To define the mirror family, we choose an
isotropic
primitive vector $f$ in the orthogonal complement $M^\perp$ of $M$ in $L$, and
consider the
lattice $\check M =
(\bfz f)_{M^\perp}^\perp/\bfz f$. Under certain arithmetic conditions on $M$,
the
lattice $\check M$ admits a primitive embedding in $M^\perp$, and we can define
the mirror moduli
space
 ${\bf K}_{\check M}$. Additional conditions on $M$ ensure that the moduli
spaces
${\bf K}_{M}$ and ${\bf K}_{\check M}$ are defined uniquely up to isomorphism,
and $\check{\check M} = M$.
Now if we have any complete family ${\cal F}$ of pseudo-ample $M$-polarized K3
surfaces its mirror
family is any complete family of pseudo-ample $\check M$-polarized K3 surfaces.

It turns out that the choice of the isotropic vector $f$ is
the analog of MS1. Property MS2 becomes

{\parindent = 1.5 cm
\item{MS2$'$} The dimension of the family $\cf$ is equal to the rank of the
Picard group
of a general member from the mirror family $\cf ^*$.
\par}

\noindent
In the case of K3 surfaces the Griffiths-Yukawa cubic becomes a quadratic map
$$S^2(H^1(V,\Theta_V))\to H^{0,2}(V)^{\otimes 2}$$ and
we have the following analog of MS3:

{\parindent = 1.5 cm
\item{MS3$'$} For some open subset $U$ of ${\cal F}$ near the boundary point
(determined by the choice of isotropic vector
$f$), for
any $V\in U$, the
Griffiths-Yukawa quadratic map can be canonically identified, after some choice
of basis
of $H^{0,2}(V)$, with the quadratic form on $\check M\otimes \bfc$.
\par}
\noindent
Note that in our case the Griffiths-Yukawa quadratic map becomes the
``constant''
quantum intersection
form and does not carry any information about rational
curves on the mirror surfaces. This agrees with the fact that
the quantum cohomology ring of a K3 surface coincides with the usual cohomology
ring {\bf [36]}.

\noindent
The next property is a very close analog of MS4:
{\parindent = 1.5 cm
\item{MS4$'$} The period map induces a holomorphic multivalued  mapping from
the subset
$U$ from MS3$'$
to the tube domain
${\rm Pic}(V')_{\bfr}+\sqrt {-1} C(X)^+\subset {\rm Pic}(V')_\bfc$ where
$C(X)^+$
 is the ample cone of
$V'\in \cf^*$. It is called the {\it mirror map}.
\par}

We also produce some computational evidence to support our point. The mirror
candidates
for a family of Calabi-Yau three-dimensional hypersurfaces in
toric Fano varieties are obtained by  Batyrev's construction from {\bf [5]}.
When this construction applies to K3 surfaces it is
``often'', but not always, gives our mirror family. This was
first noticed by Batyrev in a preprint version of {\bf [5]}. For example, in
the case when $\cf$ is the moduli family of
quartic surfaces
($M = \bfz e, (e,e) =  4$) the mirror family $\cf^*$ is the one-dimensional
family of
K3 surfaces obtained by dividing
the surfaces
$$\lambda x_0x_1x_2x_3+x_0^4+x_1^4+x_2^4+x_3^4 = 0$$
by a symplectic action of the group $(\bfz/4)^2$.
Note the analogy with Green-Plesser's first discovered example of mirror
symmetry for quintic
hypersurfaces (see {\bf [11]}, pp.1-30). In this case assertion MS3$'$ was
verified in
{\bf
[27]}.

Other examples of our computations include the mirror families for the
families of polarized K3 surfaces of degree $2n$,
for K3-covers of
Enriques surfaces and Kummer surfaces, and the fourteen families coming from
exceptional unimodal
surface singularites. In the case of polarized K3 surfaces of degree $n$ we
compute
the monodromy group of the period and mirror mappings. By other methods this
computation was done in {\bf [27]} ($n = 4$) and {\bf [21]} ($n = 2,4,6$). We
prove that the mirror moduli
space is isomorphic to the modular curve $X_0(n)^+ = H/\Gamma_0(n)^+$, where
$\Gamma_0(n)^+$ is the Fricke double extension of the modular
group $\Gamma_0(n)$. We relate the surfaces from this family to the Kummer
surfaces Kum$(E\times E')$, where
$(E,E')$ is a pair of isogeneous  elliptic curve.

There is another view of mirror symmetry more relevant to the string theory.
Here one considers the moduli
space of pairs $(V,\alpha)$ where $V$ is a Calabi-Yau manifold, and $\alpha =
B+iK\in H^2(V,\bfr/\bfz)+i{\cal K}_{V}$ with
$K$ equal to a K\"ahler form on $V$ with respect to the complex structure of
$V$. Then the mirror
mapping extends to a map of this space to itself defined on pairs $(V,\alpha)$
such that $V$ is close to
a boundary point of the moduli space of complex structures and the imaginary
part of $\alpha$ can be
represented by the Einstein-K\"ahler metric of sufficiently large volume.
The work of P. Aspinwal and D. Morrison {\bf [2,3]} treats mirror symmetry for
K3 surfaces from this point of view.

\bigskip
My acknowledgments go to V. Batyrev, D. Morrison, V. Nikulin and A. Todorov who
shared my
belief that Arnold's Strange Duality
must be the pre-history of mirror symmetry and left it up to me to
elaborate on this subject. I am also grateful to V. Batyrev, A.
Greenspoon and D. Morrison for providing some references to the relevant
literature.

\vglue .3 in

{\bf 1. Lattice polarized K3 surfaces.}
Let $X$ be a complex algebraic  K3 surface, a nonsingular projective algebraic
surface over ${\bf C}$ with
vanishing canonical class and first Betti number. Via interesection
form
the second cohomology group
$H_2(X,{\bf Z})$ is equipped with the structure of a lattice (= a free abelian
group together with an integral symmetric bilinear form on it). It is
isomorphic to the lattice $L$ equal to the orthogonal sum of three copies of
the standard hyperbolic plane $U$ (= an even unimodular indefinite lattice of
rank 2) and two copies of the lattice
$E_8$ (= an even unimodular negative definite lattice of rank 8). The
lattice $L$ is called the K3-lattice. The Picard group
${\rm Pic}(X)$ of divisor classes of $X$ is naturally identified with the
sublattice of algebraic
cycles of $H_2(X,{\bf Z})$. The unimodularity of $H_2(X,{\bf Z})$ (= Poincar\'e
duality)
allows one to identify $H_2(X,{\bf Z})$ with the second cohomology group
$H^2(X,{\bf Z}) = {\rm Hom}(H_2(X,{\bf Z}),{\bf Z})$
equipped with the lattice structure by means of the cup-product. Let
$$c:{\rm Pic}(X)\to H^2(X,{\bf Z})$$
be the corresponding injection. If one uses the interpretation of ${\rm
Pic}(X)$ as the group of isomorphism classes of
line bundles on $X$, the map $c$ corresponds to the first Chern class map. In
virtue of
the Hodge Index theorem, the lattice
${\rm Pic}(X)$ is of signature $(t_+,t_-) = (1,t)$.

The complex structure on $X$ defines the Hodge decomposition
$$H^2(X,{\bf C}) = H^{2,0}(X)\oplus H^{1,1}(X)\oplus H^{0,2}(X),$$
where
$H^{p,q}(X) \cong H^q(X,\Omega_X^p).$
It is known that complex conjugation sends $H^{2,0}(X)$ to $H^{0,2}(X)$ and
$$P_X = (H^{2,0}(X)\oplus H^{0,2}(X))\cap H^2(X,{\bf R})$$
is a positive definite 2-plane in $H^2(X,{\bf R})$. The subspace
$$H_{\bf R}^{1,1}(X) = H^{1,1}(X)\cap H^2(X,{\bf R})$$
has signature $(1,19)$.
The cone
$$V(X) = \{x\in H_\bfr^{1,1}(X)\cap H^2(X,{\bf R}): (x,x) > 0\}$$
consists of two connected components. We denote by  $V(X)^+$ the component
which
contains the
class of some K\"ahler form on $X$
with respect to the complex structure of $X$.
Let
$$\Delta(X) = \{\delta\in {\rm Pic}(X):(\delta,\delta) = -2\}.$$
By Riemann-Roch, $\Delta(X) = \Delta(X)^+\coprod \Delta(X)^-$ where
$\Delta(X)^+$ consists of
effective classes and $\Delta(X)^- = -\Delta(X)^+$. Let $W(X)$ be the subgroup
of the orthogonal group of $H^2(X,{\bf Z})$ generated by reflections
in elements from $\Delta(X)$. This group acts properly discontinuously in
$V(X)^+$ with fundamental domain
$$C(X) = \{x\in V(X)^+:(x,\delta) \ge 0, \quad\hbox{for any $\delta\in
\Delta(X)^+$}\}.$$
The set $C(X)^+$ of its interior points is the the K\"ahler cone of $X$ ({\bf
[13]}, Expos\'e X).

By the Lefschetz Theorem
$${\rm Pic}(X) = H_{\bf R}^{1,1}(X)\cap H^2(X,{\bf Z}).$$
We set
$${\rm Pic}(X)^{+} = C(X)\cap
H^2(X,{\bf Z}),\ \ \ {\rm Pic}(X)^{++} = C(X)^+\cap
H^2(X,{\bf Z}).$$
The elements of ${\rm Pic}(X)^{+}$ are
pseudo-ample  divisor classes, i.e., numerically effective
divisor classes with positive
self-intersection. Elements of ${\rm Pic}(X)^{++}$ are ample divisor classes.
Elements from $V(X)^+\cap H^2(X,{\bf Z})$ are just effective divisor
classes with positive self-intersection.

\medskip
 Now let $M$ be an even non-degenerate lattice of signature $(1,t)$.
The cone
$$V(M) = \{x\in M_{\bf R}: (x,x) > 0\} \subset M_{\bf R}$$
consists of two connected components.
We fix one of them and denote it by $V(M)^+$.

Let
$$\Delta(M) = \{\delta\in M: (\delta,\delta) = -2\}.$$
We fix a subset $\Delta(M)^+$ such that
\item{(i)}$\Delta(M) = \Delta(M)^+\coprod \Delta(M)^-$, where $\Delta(M)^- =
\{-\delta:\delta\in \Delta(M)^+\}$;
\item{(ii)} if $\delta_1,\ldots,\delta_k\in\Delta(M)^+$ and $\delta=\sum
n_i\delta_i$ with
$n_i\ge 0$ then $\delta\in \Delta(M)^+$.

The choice of subset $\Delta(M)^+$ as above defines the subset
$$C(M)^+ = \{h\in V(M)^+\cap M: (h,\delta) > 0 \quad\hbox{for all $\delta\in
\Delta(M)^+$}\}.$$

\medskip\noindent
\ \ \ {\bf Definition}. An {\it $M$-polarized}
K3 surface is a pair $(X,j)$ where $X$ is a K3 surface and $j:M\hookrightarrow
{\rm Pic}(X)$ is a
primitive lattice embedding. We say that $(X,j)$ is {\it pseudo-ample}
(resp. {\it ample}) $M$-polarized if
$$j(C(M)^+)\cap{\rm Pic}(X)^{+}\ne \emptyset$$
(resp.
$$j(C(M)^+)\cap{\rm Pic}(X)^{++}\ne \emptyset).$$

Two $M$-polarized K3 surfaces $(X,j)$ and $(X',j')$ are called isomorphic if
there
exists an isomorphism of K3 surfaces $f:X'\to X$ such that $j = f^*\circ j'$.

\medskip
{\bf Remarks (1.1).}
Note that for any pseudo-ample $M$-polarized K3 surface $(X,j)$ we have
\item{(i)} $j(V(M)^+) \subset V(X)^+$;
\item {(ii)}$j(\Delta(M)^+) = j(M)\cap \Delta(X)^+$.

Conversely, if these conditions are satisfied, and $j(M) = \pic(X)$, then
$(X,j)$ is ample
$M$-polarized.

Finally observe that a pseudo-ample $M$-polarized K3 surface is algebraic.

\medskip
{\bf Example (1.2).} Let $M = <2n>:= \bfz e, (e,e) = 2n$. Assume $n > 0$.
Choose
$V(M)^+$ to be one of the two rays in
$M_\bfr\setminus \{0\}$. A pseudo-ample \pol K3 surface $(X,j)$ is called
a {\it degree $2n$
pseudo-polarized K3 surface}. Consider the complete linear system
$|j(e)|$ and let $f$ be a rational map defined by this linear system. Then one
of the
following three possible cases occurs:
\item{(i)} (Unigonal case) $|j(e)|$ has a base curve $C\cong {\bf P}^1$,
$|j(e)-C| = |(n+1)E|$ where
$E$ is an elliptic curve. The map $f$ is a regular map from
$X$ to ${\bf P}^{n+1}$ whose image is a normal rational curve of degree $n+1$.
\item{(ii)} (Hyperelliptic case) $|j(e)|$ has no base points and $f$ is a
morphism of degree $2$ onto a normal surface of degree  $n$ in ${\bf P}^{n+1}$.
Its singular points (if any) are rational double points.

\item{(iii)}$|j(e)|$ has no base points and $f$ is a
morphism of degree $1$ onto a normal surface of degree  $2n$ in ${\bf
P}^{n+1}$.
Its singular points (if any) are rational double points.

Moreover, if $j(e)$ is ample, the unigonal case may occur only if $n = 1$. Also
in cases (ii) and
(iii), the morphism $f$ is finite and
its image is nonsingular (see {\bf [13]},
Expos\'e IV).

\vglue .3 in

{\bf 2. Local deformations.}
Fix an $M$-polarized K3 surface $(X_0,j_0)$. Let  $S$ be the local moduli space
for $X_0$. It is smooth of dimension 20 with
all Kodaira-Spencer mappings
$$\rho_s:T_s(S)\to H^1(X_s,\Theta_{X_s})$$
being isomorphisms. Let
$$H^1(X_0,\Theta_{X_0})\otimes H^1(X_0,\Omega_{X_0}^1)\to H^2(X_0,{\cal
O}_{X_0})\leqno(1)$$
be the natural pairing induced by the duality map $\Theta_{X_0}\otimes
\Omega_{X_0}^1 \to {\cal O}_{X_0}$.
As was explained in the previous section, we can view ${\rm Pic}(X)$ as a
subgroup of $H^{1,1}(X_0).$ We denote by
$$H^1(X_0,\Theta_{X_0})_{j_0}$$
the orthogonal complement of $j_0(M)$ with respect to the
pairing $(1)$. Let
$$H^1(X_0,\Omega_{X_0}^1)_{j_0} = H^1(X_0,\Omega_{X_0}^1)/j_0(M).$$
In view of the canonical pairing
$$H^{1,1}(X_0)\otimes H^{1,1}(X_0)\to H^{2,2}(X_0) \cong {\bf C}$$
we may identify $H^1(X_0,\Omega_{X_0}^1)_{j_0}$ with
$$H^{1,1}(X_0)_{j_0} : = (j_0(M)_{\bf C})^\perp_{H^{1,1}(X_0)}.$$
The pairing $(1)$ induces the map
$$d_1:H^1(X_0,\Theta_{X_0})_{j_0}\to {\rm
Hom}(H^1(X_0,\Omega_{X_0}^1)_{j_0},H^2(X_0,{\cal O}_{X_0})).\leqno(1')$$

\medskip

\proclaim \ \ \ Proposition (2.1). There is a local moduli space $S_M$ of
isomorphism
classes of
$M$-polarized K3 surfaces. It is smooth of dimension $19-t$.  Its tangent space
at each point
$s\in S$ is naturally isomorphic to $H^1(X_0,\Theta_{X_0})_{j_0}$.

{\sl Proof}. In the case $t = 0$ this is a theorem from ${\bf [15]}$. The
general
case is proved similarly.

\medskip
Let
$$H^1(X_0,\Theta_{X_0})\otimes H^0(X_0,\Omega^2_{X_0})\to
H^1(X_0,\Omega^1_{X_0})\leqno(2)$$
be the natural pairing induced by the contraction map $\Theta_{X_0}\otimes
\Omega^2_{X_0}\to \Omega^1_{X_0}$.
Composing  $(2)$ with the projection $H^1(X_0,\Omega^1_{X_0}) \to
H^1(X_0,\Omega^1_{X_0})_{j_0}$ and
restricting the composition to $H^1(X_0,\Theta_{X_0})_{j_0}$ we get the map
$$d_2:H^1(X_0,\Theta_{X_0})_{j_0}\to {\rm
Hom}(H^0(X_0,\Omega^2_{X_0}),H^1(X_0,\Omega^1_{X_0})_{j_0}).\leqno(2')$$
Let
$$(d_1,d_2):H^1(X_0,\Theta_{X_0})_{j_0}\to {\rm
Hom}(H^{1,1}(X_0)_{j_0},H^{0,2}(X_0))\oplus {\rm
Hom}(H^{2,0}(X_0),H^{1,1}(X_0)_{j_0}).$$
The formula
$$(\theta_1,\theta_2) = d_1(\theta_2)\circ d_2(\theta_1):H^{2,0}(X_0)\to
H^{0,2}(X_0) $$
defines the bilinear form
$$H^1(X_0,\Theta_{X_0})_{j_0}^{\otimes 2}\to {\rm
Hom}(H^{2,0}(X_0),H^{0,2}(X_0)).$$
The canonical pairing
$$H^{2,0}(X_0)\otimes H^{0,2}(X_0) \to H^{2,2}(X_0) = H^4(X_0,{\bf C}) \cong
{\bf C}$$
allows one to identify the space of values of the bilinear form with the space
$H^{0,2}(X_0)^{\otimes 2}$.
 One can check that this pairing is symmetric {\bf [16]} and hence defines
the linear map
$${\rm Yu}: S^2(H^1(X_0,\Theta_{X_0})_{j_0})\to H^{0,2}(X_0)^{\otimes 2}$$
which we call the {\it Griffiths-Yukawa quadratic map} for $M$-polarized
K3 surfaces. A choice of an
isomorphism
$$H^{0,2}(X_0)\cong {\bf C}$$
allows one to identify the map Yu with a quadratic
form on the space $H^1(X_0,\Theta_{X_0})_{j_0}$.

The Griffiths-Yukawa quadratic map can also be expressed in terms of the
intersection form on $H^{1,1}(X_0)$ as follows.
First observe that the map
$$d_1:H^1(X_0,\Theta_{X_0})_{j_0} \to {\rm
Hom}(H^{2,0}(X_0),H^{1,1}(X_0)_{j_0}) \cong H^{0,2}(X_0)\otimes
H^{1,1}(X_0)_{j_0}\leqno (3)$$
is bijective (since it is injective and both spaces have the same dimension).
The pairing
$$ H^{1,1}(X_0)_{j_0}\otimes H^{1,1}(X_0)_{j_0} \to {\bf C}$$
defines the symmetric bilinear map
$$(H^{2,0}(X_0)^*\otimes H^{1,1}(X_0)_{j_0})\otimes (H^{2,0}(X_0)^*\otimes
H^{1,1}(X_0)_{j_0})\to H^{0,2}(X_0)^{\otimes 2} .\leqno (4)$$
Then it follows from the definition of the Griffiths-Yukawa quadratic map that
$${\rm Yu} = {\rm Yu}'\circ d_1,$$
where
$${\rm Yu}': H^{2,0}(X_0)\otimes H^{1,1}(X_0)_{j_0} \to H^{0,2}(X_0)^{\otimes
2}$$
is the quadratic map defined by $(4)$.

\vglue .3 in

{\bf 3. The period map.}
The map $(3)$  can be naturally interpreted as the differential of
the period mapping for \pol K3 surfaces. Let $M$ be a lattice of signature
$(1,t)$ embeddable into the K3-lattice $L$.
Fix a lattice embedding $i_M:M \to L$. We shall often identify $M$ with the
image $i_M(M)$.   Let
$$N =M^\perp_L $$
be the orthogonal complement of $M$ in $L$. It is a lattice of signature
$(2,19-t)$.

\medskip\noindent
\ \ \ {\bf Definition}. A {\it marked \pol K3 surface} is a pair $(X,\phi)$,
where $X$
is a K3 surface
together with an isomorphism of lattices
$\phi:H^2(X,{\bf Z})\to L$ such that $\phi^{-1}(M)\subset {\rm Pic}(X)$.
The
pair  $(X,j_\phi)$ with
$j_\phi = \phi^{-1}|M: M\to {\rm Pic}(X)$ is a $M$-polarized K3 surface.
In this way we can speak about pseudo-ample and ample
marked \pol K3 surfaces. Two marked surfaces $(X,\phi)$ and $(X',\phi')$ are
called isomorphic if there
exists an isomorphism of surfaces $f:X\to X'$ such that $\phi' = \phi\circ
f^*$.

\medskip
Given a marked \pol K3 surface $(X,\phi)$, the Hodge decomposition of
$H^2(X,{\bf C})$ defines
the point $\phi(H^{2,0}(X))$ in ${\bf P}(L_{\bf C})$. Since $H^{2,0}(X)$ is
orthogonal to
$H^{1,1}(X)$ (with respect to the cup-product in $H^2$), the line
$\phi(H^{2,0})$ is always orthogonal to
$\phi(j(M)) = M$. Therefore
$$\phi(H^{2,0}(X))\in {\bf P}(N_{\bf C}) \subset {\bf P}(L_{\bf C}).$$
Let $Q$ be the quadric in ${\bf P}(N_{\bf C})$ corresponding to the quadratic
form
on $N_{\bf C}$ defined  by the lattice $N$. For any $\omega\in H^{2,0}(X)$ we
have
$$(\omega,\omega) \in H^{4,0}(X) = \{0\}.$$
This shows that $\phi(H^{2,0})\in Q$. Also
$$(\omega,\bar \omega)\in \bfr_+ \subset H^{2,2}(X)\cap H^4(X,\bfr) = \bfr.$$
Therefore  $\phi(H^{2,0}(X))$ is contained in an open (in the usual topology)
subset $D_M$ of the quadric $Q$
defined by the inequality $(\omega,\bar \omega) > 0$.
By assigning to $H^{2,0}(X_0)$ the positive definite 2-plane
$P_X \subset N_\bfr$ together with the orientation defined by the choice of the
isotropic line $H^{2,0}(X_0)\subset P_X\otimes \bfc$,
we can identify $D_M$ with the symmetric homogeneous space
$O(2,19-t)/SO(2)\times O(19-t)$ of
oriented positive definite 2-planes in $N_\bfr$. The space consists of two
connected components
each isomorphic to a bounded Hermitian domain of type $IV_{19-t}$. The
involution which
switches the two components is induced by the complex conjugation map $Q\to Q$.
We shall call the point $\phi(H^{2,0})\in D_M$ the {\it period point} of
$(X,\phi)$.

\medskip
Let $S_M$ be the local moduli space of an \pol K3 surface $(X_0,j_0)$. Since
$S_M$ is contractible, we can
choose a marking $\phi:H^2(X_s,{\bf C})\to L$ for all $X_s, s\in S_M$. We fix
this marking and identify
$H^2(X_s,{\bf C})$ with $L_{\bf C}$. The complex structure on $X_s$ defines
the Hodge decomposition and hence the period point $H^{2,0}(X_s)\in D_M$.
By the  Local Torelli Theorem ({\bf [13]}, Expos\'e V) the resulting period map
$$p: S_M \to D_M$$
is a holomorphic map which is locally an isomorphism in a neighborhood of the
origin
$0 = (X_0,j_0)$. Let $\mu = p(0)\in D_M$ be the period point of $(X_0,j_0)$. We
shall identify it
with the subspace $H^{2,0}$ of $N_\bfc$. Then
$$T_\mu(D_M) \cong Hom(\mu,\mu^\perp/\mu) =
Hom(H^{2,0}(X_0),H^{1,1}(X_0)_{j_0})\cong  H^{0,2}(X_0)\otimes
H^{1,1}(X_0)_{j_0}.$$
The differential of the period map
$$dp_0:T_0(S_M)\to T_\mu(D_M)$$
is the bijective map
$d_1$ from $(3)$.

\medskip
Let ${\cal K}_M$ be the fine moduli space of marked \pol K3 surfaces. It is
obtained by
gluing local moduli spaces of marked \pol K3 surfaces and is a (non-separated)
analytic space
(see {\bf [13]}, Expos\'e XIII, {\bf [30]}). The local period maps are glued
together to give a holomorphic map
$$p:{\cal K}_M\to D_M.$$
This map is the restriction of the period map $P:{\cal M}\to D \supset D_M$ for
marked K\"ahler K3 surfaces.
According to the Global Torelli Theorem of Burns-Rappoport and the Surjectivity
Theorem of Todorov
the latter map is \'etale and surjective (see loc.cit). The former theorem also
describes
the fibres of the period map.

Let $(X,\phi)$ be a marked \pol K3 surface. Then the image of the data
$$(P_X,V(X)^+,\Delta(X)^+,C(X))$$
under the map $\phi$ defines the data
$(\pi, V_\pi^+,\Delta_\pi^+, C_\pi),$
where
\item{(i)} $\pi$ is a positive oriented 2-plane in $N_\bfr$;
\item{(ii)}$V_\pi^+$ is a connected component of the cone $\{x\in
\pi^\perp:(x,x) > 0\}$
;
\item{(iii)}$\Delta_\pi^+$ is a subset of
$\Delta_\pi = \{e\in \pi^\perp\cap L: (e,e) = -2\}$ such that
$\Delta_\pi = \Delta_\pi^+\coprod -\Delta_\pi^+$;
\item{(iv)}$C_\pi = \{x\in V_\pi^+ :(x,e) \ge 0$ for any $e\in
\Delta_\pi^+\}.$

Note that $V_\pi^+$ is uniquely determined by $\pi$ (since $V(M)^+$ is fixed)
and $C_\pi$ is
determined by $\Delta_\pi^+$.

\medskip

\proclaim\ \ \ Theorem (3.1). The restriction of the  period map $p:{\cal K}_M
\to D_M$ to the subset
${\cal K}_M^{pa}$ of isomorphism classes of marked pseudo-ample $M$-polarized
K3 surfaces
is surjective. There is a natural bijection between the fibre
of the map
$$p':{\cal K}_M^{pa}\to D_M$$
over a point $\pi\in D_M$ and the subgroup $W_\pi(N)$ of
isometries of $L$ generated
by reflections in vectors from $\Delta_\pi\cap N$.

{\sl Proof}. It follows from the Global Torelli Theorem that the fibre
$p^{-1}(\pi)$ is bijective to the set of possible pairs
$(V_\pi^+,\Delta_\pi^+)$. The group
$W_\pi(L)\times \{\pm 1\}$, where $W_\pi(L)$ is generated by reflections in
elements from $\Delta_\pi$, acts
transitively on the
fibre.
Pick up a point $(X,\phi)$ in $p^{-1}(\pi)$ corresponding to
$(V_\pi^+,\Delta_\pi^+)$ such
that $V(M)^+\subset V_\pi^+$. Let $h\in C(M)^+$, we may choose $(X,\phi)$ such
that
$h\in C_\pi$. This is possible
because $C_\pi$ is a fundamental domain for the action of $W_\pi(L)$ in
$V_\pi^+$. Since $\pi\in D_M$,
$\phi^{-1}(M)\subset \pic(X)$ and $j_\phi(h)\in \pic(X)^+$. Composing $\phi$
with some reflections
from $\Delta(M)$, we may assume that $h\in j_\phi(C(M)^+)$. Thus $(X,\phi)$ is
a marked pseudo-ample $M$-polarized K3 surface with $p((X,\phi)) =\pi$.
This proves the surjectivity.

Let $(X,\phi)\in {\cal K}_M^{pa}$ then the fibre
of
$p':{\cal K}_M^{pm} \to D_M$ over
$\pi = p((X,\phi))$ is bijective to the set of subsets $\Delta_\pi^+$ such that
$\Delta_\pi^+\cap M = C(M)^+$. The stabilizer $G$ of this set in $W_\pi(L)$ is
the subgroup
$W_\pi(N)$. In fact, it follows from the properties of reflection groups that
$G$ is generated by
reflections in vectors $\delta\in \Delta_\pi$ such that $(\delta,h) = 0$, where
$h\in C(M)^+$.
Since $C(M)^+$ linearly spans $M_{\bf R}$ we must have $\delta\in N$. This
proves the theorem.

\medskip
For any
$\delta\in \Delta(N) = \{x\in N: (x,x) = -2\}$, set
$$H_\delta = \{z\in N_{\bf C}: (z,\delta) = 0\},$$
$$D_M^\circ = D_M\setminus(\bigcup_{\delta\in \Delta(N)}H_\delta\cap D_M).$$
Let $(X,\phi)$ be an ample marked \pol K3 surface. Then $(j_\phi(M))^\perp\cap
H^{1,1}$ cannot contain vectors $v$ with
$(v,v) = -2$. This shows that the period point $\pi = \phi(P_X)$
satisfies
$\Delta_\pi\cap N = \emptyset$. This implies the following:

\medskip
\proclaim \ \ \ Corollary (3.2). Let ${\cal K}^a_M$ denote the subset of ${\cal
K}_M$
which consists of isomorphism classes of
marked ample \pol K3 surfaces. The period map induces a bijective map
$$p:{\cal K}^a_M \tilde\to D_M^\circ.$$

\medskip
Next we want to get rid of markings of \pol K3 surfaces. For any lattice $S$ we
denote by $O(S)$ its orthogonal group.
The group
$$\Gamma(M) = \{\sigma\in O(L): \sigma (m) = m \quad\hbox{for any $m\in
M$}\}.$$
acts on the moduli space ${\cal K}_M$ transforming $(X,\phi)$ to
$(X,\phi\circ\sigma)$
without changing
the isomorphism class of the \pol K3 surface $(X,j_\phi)$.

Let $\Gamma_M$ be the image of $\Gamma(M)$ under the natural injective
homomorphism
$$\Gamma(M) \to O(N).$$

\smallskip
\proclaim \ \ \ Proposition (3.3). Let $A(N) = N^*/N$ be the discriminant group
of the lattice $N$,
and let
$O(N)\to {\rm Aut}(A(N))$ be the natural homomorphism. Denote its kernel by
$O(N)^*$.
Then
$$\Gamma_M = O(N)^*.$$
In particular, $\Gamma_M$ is a subgroup of finite index
in $O(N)$.

{\sl Proof}. This follows from Corollary 1.5.2 in {\bf [31]}.
\medskip

The group $O(N)$ is an arithmetic subgroup of $O(2,19-t)$ and by the previous
proposition so is $\Gamma_M$.
Since $D_M$ is the union of two copies of a Hermitian homogeneous domain, the
quotient
$D_M/\Gamma_M$ is a quasi-projective algebraic variety with at most two
irreducible
components.  Obviously $\Gamma(M)$ contains the subgroup of
$O(L)$
generated by reflections in
vectors $\delta\in N$ with $(\delta,\delta) = -2$. Thus each fibre of the map
$p': {\cal K}_M^{pa}\to D_M$ is mapped to the same orbit in ${\bf
K}_M/\Gamma(M)$. Applying Theorem 3.1, we
obtain that
the period mapping descends to a bijection
$${\cal K}_M^{pa}/\Gamma(M) \cong  D_M/\Gamma_M.$$
Since the elements of the quotient set ${\cal K}_M^{pa}/\Gamma(M)$ are
isomorphism classes of
pseudo-ample $M$-polarized K3 surfaces, we are able to endow the set ${\cal
K}_M^{pa}/\Gamma(M)$ with a
structure of a quasi-projective algebraic variety. We denote this variety by
${\bf K}_M$.

Assume that $M$ satisfies the following condition:

\item {(U)} For any two primitive embeddings $i_1,i_2 :M\hookrightarrow L$,
there exists an isometry
$\sigma:L\to L$ such that $i_1\circ\sigma = i_2$.

Let $(X,j)$ be a pseudo-ample $M$-polarized  K3 surface. Take any marking
$\phi:H^2(X,\bfz)$
$\to L$. Composing it with
$j:M\to H^2(X,\bfz)$, we obtain a primitive embedding $i':M\to L$. Replacing
$\phi$ with $\phi\circ \sigma$ for appropriate
isometry $\sigma$ of $L$, we obtain a new marking $\phi':H^2(X,\bfz)\to L$ such
that $j= j_{\phi'}$. This
shows that any isomorphism class of $(X,j)$ is represented by a point of ${\bf
K}_M$. So we may view ${\bf K}_M$ as the moduli space
of $M$-polarized K3 surfaces.

Similarly we can define the variety
${\bf K}_M^{a}$ of isomorphism classes of ample $M$-polarized K3 surfaces. We
have
$${\bf K}_M^a =  D_M^\circ/\Gamma_M.$$
Observe that $O(N)$ has only finitely many orbits in the set of primitive
vectors with given value
of the quadratic form (for example, this follows from Proposition 1.15.1 of
{\bf [31]}). This shows that the complement of
$D_M^\circ/\Gamma_M$ in $D_M/\Gamma_M$ is the union of finitely many
hypersurfaces, in particular,
$D_M^\circ/\Gamma_M$ is an open Zariski subset of $D(M)/\Gamma_M$.

\medskip
{\bf Remark (3.4).} In fact, following {\bf [13]}, \'expose XIII, one can show
that $D_M/\Gamma_M$
is a coarse moduli space of pseudo-ample $M$-polarized K3 surfaces.
First we define a family of \pol K3 surfaces. This is a family $f:{\cal X}\to
S$ of K3 surfaces
together with a homomorphism of
sheaves $M_S\to {\cal P}ic_{{\cal X}/S}$
where ${\cal P}ic_{{\cal X}/S}\subset R^2f_*(\bfz)$ is the relative Picard
sheaf. We can define a family of
pseudo-ample \pol K3 surfaces by requiring
additionally that each $(f^{-1}(s),j_s)$ is pseudo-ample \pol K3 surface.
Since
${\cal K}_M$ is a fine moduli space for marked \pol K3 surfaces, a family
$(f:{\cal X}\to S,\phi)$ of marked \pol K3 surfaces is equivalent to a
holomorphic map $\alpha:S\to {\cal K}_M$. Composing this map with the
period we obtain that  $(f:{\cal X}\to S,\phi)$ defines a holomorphic map $\bar
\alpha:S\to D_M/\Gamma_M$.
Given a family $(f:{\cal X}\to S,j)$ of pseudo-ample \pol K3 surfaces, after
localizing $S$, we equip it with marking, and define the map
$\bar \alpha:S\to D_M/\Gamma_M$ which does not depend on the choice of the
marking. When $S$ is a point, we get a
a bijection ${\cal K}_M^{pa}/\Gamma(M) \cong  D_M/\Gamma_M$. This proves that
$D_M/\Gamma_M$ is
a coarse moduli space.
Similarly we prove that
$D_M^\circ/\Gamma_M$ is a coarse moduli space for ample \pol K3 surfaces.

I don't know any algebraic construction for ${\bf K}_M$ except when $M$ is of
rank 1.

\vglue .3 in
{\bf 4. Tube domain realization of the period space.}
Let $b:V\times V\to {\bf C} $ be a non-degenerate symmetric bilinear form on a
complex vector
space and let $Q: b(x,x) = 0$ be the corresponding non-degenerate quadric in
the projective
space ${\bf P}(V)$ associated to $V$. For any non-zero vector $v\in V$ we
denote by $[v]$ the line
${\bf C}v \in {\bf P}(V)$. For any $v\in V\setminus \{0\}$ the hyperplane
$H_v = \{w\in V: b(w,v) = 0\}$ intersects $Q$ along the quadric
$$Q(v) = Q\cap H_v = \{x\in Q: v\in PT_x(Q)\} \subset H_v,$$
where $PT_x(Q)$ is the projective tangent space of $Q$ at the point $[v]$. If
$b(v,v) = 0$, i.e.
$[v]\in Q$, the hyperplane $H_v$ coincides with $PT_{[v]}(Q)$  and $Q(v)$ is
the cone over
the quadric $\bar Q(v) \subset {\bf P}(H_v/{\bf C}v)$ with the vertex at $[v]$.
In other words,
the projection map $Q\setminus \{[v]\}\to {\bf P}(V/{\bf C}v)$ is an
isomorphism outside $Q(v)$, and blows down
$Q(v)\setminus\{[v]\}$ to the quadric $\bar Q(v)$.

\medskip
We shall apply the previous remarks to our situation where
$V = N_{\bf C} \subset L_{\bf C}$ with the symmetric bilinear form defined by
the lattice
$N = M^\perp$. The period space $D_M$ is a subset of the quadric $Q$ defined by
the inequality
$(\mu,\bar \mu) > 0$.

\medskip
\proclaim \ \ \ Lemma (4.1).  Let $f\in (M^\perp)_\bfr, (f,f) = 0$. Then
$$D_M \cap Q([f]) = \emptyset.$$

{\sl Proof.} Suppose there exists $\mu\in D_M\cap Q([f])$.
 Since $f$ is a real vector, and the bilinear form originates from the
lattice structure, we have  $(\bar \mu,f) = (\mu,f) = 0$. This implies
that
$f\in P^\perp$, where
$P \subset N_{\bf R}$ is the positive definite 2-plane spanned by the real and
imaginary part of $\mu$. However,
the signature $(t_+,t_-)$ of the space $(M^\perp)_{\bf R}$ satisfies $t_+ = 2$.
Therefore
$P^\perp$ is negative definite and does not contain isotropic vectors. This
contradiction proves the assertion.

\medskip
 From now on we assume that $t\le 18$, i.e., rank$(M) \le 19.$ This ensures
that the lattice $N = M^\perp$ is indefinite.
Let us fix an isotropic vector $f\in N_\bfr$. We set
$$W_f =  \{x\in N_\bfr: (x,f) = 1\}/\bfr f,$$
$$V_f = \{x\in N_\bfr: (x,f) = 0\}/\bfr f.$$
By Lemma (4.1), the projection map $\pi:Q\setminus Q([f]) \to  \bfp(N_\bfc/{\bf
C}f)$ maps $D_M$ isomorphically
onto a subset of the affine space
$$A_f = \bfp(N_\bfc/{\bf C}f)\setminus {\bf P}((V_f)_\bfc) \cong \{z\in N_\bfc:
(z,f) = 1\}/\bfc f =$$
$$ W_f+iV_f = \{z = x+iy\in N_\bfc/\bfc f = (N_\bfr/\bfr f)+i(N_\bfr/\bfr f):
(x,f) = 1, (y,f) = 0\}.$$

\medskip
\proclaim \ \ \ Theorem (4.2).
The projection map
$Q\setminus Q([f]) \to {\bf P}(N_\bfc/{\bf C}f)$ defines an analytic
isomorphism
$$D_M \cong \{x+iy\in A_f: (y,y) > 0\}.$$

{\sl Proof}. This is just the translation of the condition $(\mu,\bar
\mu) > 0$ in terms
of the projection map. We write any $\mu\in D_M$ in the form
$\mu = \lambda f+x+iy$ where $x +\bfr f\in W_f, y+\bfr f\in V_f$. We have
$$0 = (\mu,\mu) =  [2\Re(\lambda)+(x,x)-(y,y)]+i[2\Im(\lambda)+2(x,y)],$$
hence $(x,x)-(y,y)+2\Re(\lambda) = 0$. This implies that
$$0 < (\mu,\bar \mu) =  2\Re(\lambda)+(x,x)+(y,y) = 2(y,y).$$
This proves the assertion.

\medskip
Recall that for any real affine space $W$ with the translation space $V$ and an
open connected
cone $C\subset V$ which does not contain an affine line, the set
$$\Omega(W,V,C) = \{z = x+iy: x\in W, y\in C\} \subset W_\bfc$$
is called the {\it tube domain} associated to the cone $C$ in $V$.
In the special case when $V$ is equipped with a nondegenerate quadratic form
with signature $(1,n)$ and the cone
$V^+$  is one of the two connected components of the cone $\{x\in V: (x,x) >
0\}$ the tube
domain $\Omega(W,V,C)$ is a bounded Hermitian symmetric domain of type
$IV_n$. This can be applied to our situation where
$V = V_f$. Fix a connected component $V_f^+$ of the cone $\{x\in V_f:(x,x)>
0\}$. Restricting the period map to a connected component $D_M^+$ of $D_M$, we
obtain

\medskip
\proclaim \ \ \ Corollary (4.3). The choice of an isotropic vector $f\in
N_\bfr$
defines an
isomorphism
$$D_M^+ \cong \Omega(W_f,V_f,V_f^+).$$

\medskip
\proclaim \ \ \ Corollary (4.4). For any $\mu\in D_M$ the choice of a
representative
$\ell\in L$ of $\mu$ with $(\ell,f) = 1$
defines  a canonical isomorphism
$$\alpha_\mu:T_\mu(D_M) \to (V_f)_\bfc.$$
If $[(X,\phi)]\in {\cal K}_M$ is the isomorphism class of a marked \pol
K3 surface with the period point $\mu$, then
the pre-image of the quadratic form on $(V_f)_\bfc$ under the map
$$\alpha_\mu\circ dp_{[(X,\phi)]}:H^1(X,\Theta_X) \to (V_f)_\bfc$$
coincides with the Griffiths-Yukawa quadratic form with respect to the
normalization
$H^{0,2} \tilde \to \bfc$
defined by the linear function  $\phi^{-1}(\ell)\in H^{2,0}$.

{\sl Proof.} The map
$$\alpha_\mu:T_\mu(D_M) = Hom(\mu,\mu^\perp/\mu) \to (V_f)_\bfc$$
is the composition of the differential of the map $D_M\to A_f$ at the point
$\mu$ and the
differential of the translation map
$A_f \to (V_f)_\bfc, z\mapsto z-\mu$. Explicitly it sends
$\psi:\mu\to \mu^\perp/\mu$ to $\psi(\ell)'-(\psi(\ell)',f)\ell \ {\rm mod}\
\bfc f$, where $\psi(\ell)'$ is a representative of
$\psi(\ell)$ in $\mu^\perp$.

\medskip
{\bf Remark (4.5).} In general there is no canonical trivialization of the
affine
space $W_f$. However,
a choice of an isotropic vector $g\in N_{\bf R}$ with $(f,g) = 1$ defines the
trivialization
$$W_f \to (V_f), x\to x-g.$$
If we choose to identify $V_f$ with $(\bfr f+\bfr g)_{N_\bfr}^\perp$, then the
explicit
isomorphism
$\alpha:\Omega_f \to D_M$ is given by the formula:
$$\alpha(z) = {\bf C}(-{1\over 2}(z,z)f+g+z).$$

\vglue .3 in

{\bf 5. Some  arithmetical conditions on $M$.} We are going to put some
arithmetical conditions on our
lattice $M$ to ensure, for example, condition (U) in section (3).

 For each non-degenerate even lattice $S$ we
denote by $A(S) = S^*/S$ the discriminant group
of $S$ equipped with the quadratic map
$$q_S:A(S) \to {\bf Q}/2\bfz, \ \ q_M(x+S) = (x,x)+2\bfz,$$
where the bilinear form of $S$ is extended to a ${\bf Q}$-valued bilinear form
on $S^*$.

For example, for any integer $m \ne 0$, let $U(m)$ denote the lattice of rank 2
with a basis
$(e,e')$ such that
$(e,e') = m, (e,e) = (e',e') = 0$. Then $A(U(m)) = (\bfz/m\bfz)^2$ with
$q_{U(m)}$ defined by the formula:
$$q((a + m\bfz,b+m\bfz)) = {2ab\over m}+2\bfz.$$

\medskip

It is clear that for any isometry $\sigma\in O(L)$, we have a canonical
isomorphism $D_M\to D_{\sigma(M)}$ which
defines a canonical isomorphism of the moduli spaces
$${\cal K}_M \cong {\cal K}_{\sigma(M)},$$
where we choose $(V(M)^+,\Delta(\sigma(M))^+)$ to be equal to
$(\sigma(V(M)^+),\sigma(\Delta(M)^+)$.
The next result of Nikulin gives a condition implying that any two primitive
embeddings
$i:M\to L,i': M\to L$ differ by an isometry of $L$.

\medskip
\proclaim \ \ \ Proposition (5.1). Let $S$ be an even lattice of signature
$(1,t)$ with
$t\le 19$.
Assume that $l(A(S)) \le 20-t$ or $t\le 10$. Then there exists a primitive
embedding
$S\hookrightarrow L$. Moreover, such an embedding is unique up to an isometry
of $L$ if for
each prime $p\ne 2$ the
$p$-primary component $A(S)_p$ of $A(S)$ satisfies
$l(A(S)_p) \le 19-t$ and, if $l(A(S)_2 = 21-t$, $A(S)_2$ contains as a
direct summand the
discriminant form of the lattice $U(2)$.

{\sl Proof.} See {\bf [31]}, Corollary 1.12.3, Theorems 1.12.4, 1.14.4.

\medskip
\proclaim \ \ \ Corollary (5.2). Any even lattice $M$ of signature $(1,t)$ with
$t\le 9$ admits a
unique primitive embedding in the K3-lattice $L$. In particular, the moduli
space ${\bf K}_M^{a}$ of ample \pol K3 surfaces
is not empty (and of dimension $19-t$).

{\sl Proof.} In fact, ${\bf K}_M^{a}$ is a Zarisli-open non-empty subset in
$(19-t)$-dimensional algebraic variety $D_M/\Gamma_M$.

\medskip
Next we want to study primitive isotropic vectors $f$ in a non-degenerate even
lattice $S$.
Consider $f$ as an element of $S^*$ and let ${\rm div}(f)$ be the positive
generator of the image
of the linear map
$f:S\to \bfz$. Let $f^* = {1\over {\rm div}(f)}f+S \in S^*$. Clearly $f^*$ is
an
isotropic element of the discriminant quadratic form $A(S)$.
Let $I(S)$ denote the set of primitive isotropic vectors in $S$, and $I(A(S))$
be the
same for
$A(S)$.  The map $f\to f^*+S$ defines a map
$I(S)\to I(A(S))$. The orthogonal group $O(S)$ acts naturally on the source and
the target of this map, and the map
is compatible with this action.  Let
$$O(S)^* = Ker(O(S)\to O(A(S)).$$
This group acts on the fibres of the map $I(S)\to I(A(S))$.

\medskip
\proclaim \ \ \ Proposition (5.3). The map
$$I(S)/O(S)^* \to I(A(S)), f\mapsto f^* +S$$
is surjective if $S$ admits the lattice $U = U(1)$ as an orthogonal summand.
The map is bijective if
$S$ admits the lattice $U\perp U$ as an orthogonal summand.

{\sl Proof.} See {\bf [37]}, Lemmas 4.1.1 and 4.1.2.

\medskip
{\bf Definition}. An isotropic vector $f\in I(S)$ is called $m$-{\it
admissible} if
${\rm div}(f) = m$ and there exists $g\in I(S)$ with $(f,g) = m, {\rm div} g =
m$.

\medskip
\proclaim \ \ \ Lemma (5.4). The following conditions are equivalent:
\item{(i)} $f\in I(S)$ is $m$-admissible;
\item{(ii)} there exists a
primitive lattice embedding
$i:U(m)\to S$ such that $S= i(U(m))\oplus i(U(m))^\perp$ and $f\in i(U(M))$.

{\sl Proof.} $(i)\Rightarrow  (ii)$. Let $g\in I(S)$ such that $(f,g) = m$.
The sublattice $U'$ spanned by
$f$ and $g$ is primitive, contains $f$  and is isomorphic to $U(M)$. Since for
any $s\in S$,
$m$ divides $(s,f)$ and $(s,g)$, we obtain
$s-{(s,g)\over m}f-{(s,f)\over m}g\in U_S^\perp.$ This shows that $S= U'\perp
U^\perp.$

$(i)\Leftarrow  (ii)$ Obvious.

\medskip
\proclaim \ \ \ Proposition (5.5). Let $S$ be an even indefinite non-degenerate
lattice
of signature $(t_+,t-)$.
Then $S$ admits the lattice $U(m)$ as an orthogonal summand if the following
conditions are satisfied:
\item{(i)} $A(U(m))$ is isomorphic to an orthogonal summand of $A(S)$ with
respect to the bilinear
form defined by $q_S$;
\item{(ii)} $l(A(S)) \le t_++t_--3$.

{\sl Proof.} Let  $A'$ be the orthogonal complement of $A(U(m))$ in $A(S)$.
Then
$l(A') \le l(A(S)) \le {\rm rank}\ S-3 < (t_+-1)+(t_--1)$. By Corollary 1.10.2
from {\bf
[31]},
there exists a lattice $S'$ with signature $(t_+-1,t_--1)$ and $A(S') \cong
A'$.  Thus the lattice
$U(m)\perp S'$ has the same signature and the same discriminant quadratic form
as the lattice $S$.
By Corollary 1.13.3 from loc. cit. we obtain $S\cong U(m)\perp S'$.

\medskip
\proclaim \ \ \ Proposition (5.6). Suppose $M^\perp$ contains an $m$-admissible
isotropic vector with
$m\le 2$. Then the moduli space
${\bf K}_M$ is irreducible.

{\sl Proof.} Let $M^\perp = U(m)\perp  M'$. The isometry
$-{\rm id}_{U(m)}
\oplus{\rm id}_{M'}$ of $M^\perp$ acts identically on the
discriminant group of $M^\perp$, hence extends to an
isometry of $\sigma\in \Gamma(M)$ of $
L$ (see Proposition (3.3)). Obviously $\sigma$
switches the
orientation of a positive definite 2-plane $\pi \subset (M^\perp)_\bfr$ spanned
by a vector $x\in U(m)$ with $(x,x) > 0$ and a vector $y\in M'$ with $(y,y) >
0$. Hence
$\sigma$ switches the two connected
components of $D_M$. This implies that $D_M/\Gamma(M) = D_M/\Gamma_M$ is
irreducible.

\vglue .3 in

{\bf 6. Mirror symmetry.} Now we  are ready to define the mirror family.
Pick up an $m$-admissible isotropic vector $f$ in  $N=M^\perp$
(Proposition (5.3) gives some sufficient conditions for its existence).
Then  $M^\perp = U'\perp \check M$, where $U'\cong U(m)$ and $f\in U'$.
The sublattice $\check M$ is of signature $(1,18-t)$.  We
have
$$(\bfz f)^\perp_{M^\perp}/\bfz f \cong \check M.$$
So as an abstract lattice, $\check M$ does not depend on the choice of $U'\cong
U(m)$ containing
$f$.
Let us fix $(V(\check M)^+, \Delta(\check M)^+,C(\check M)^+)$ and use the
embedding $i:\check M \hookrightarrow M^\perp \subset L$
to introduce the moduli space ${\bf K}_{\check M}$ of
$\check M$-polarized K3 surfaces.

\medskip
{\bf Definition.}
The moduli space ${\bf K}_{\check M}$ is called the {\it
mirror moduli space} of
${\bf K}_{\check M}$.

\medskip
The definition depends obviously on the choice of $U'$ which determines the
embedding $\check M \hookrightarrow L$.
If we replace
$i$ with composition $i'= \sigma\circ i$ where $\sigma \in O(L)$, then
$i'(M)\subset \sigma(M^\perp) = \sigma(M)^\perp$. Thus the new ${\bf
K}_{i'(\check M)}$ will be
equal to the mirror of
${\bf K}_{\sigma(M)}\cong {\bf K}_{M}$. Thus, if we put conditions on $M$ which
ensure the
uniqueness of primitive embedding
of $\check M$ in $L$, we obtain that the isomorphism class of the mirror moduli
space
depends only on the choice of $f\in M^\perp$. Since the signature and the
discriminant group of $\check M$ can be read off
from the signature and discriminant of $M$, we can apply Proposition (5.1) to
get some sufficient
conditions on $M$ which guarantee that our construction is well-defined.

\medskip
Note the obvious relations
$${\rm dim}\  {\bf K}_{\check M}+{\rm dim}\  {\bf K}_M = 20,$$
$$ {\rm dim}\  {\bf K}_{\check M} = {\rm rank}\  M = {\rm rank\ Pic}(X),$$
where for any marking $\phi$ of $X\in {\bf K}_M$ the period of $(X,\phi)$ does
not belong to a
subvariety of the form $D_{M'}$ for some sublattice $M'$ of $L$ with
$M \subset M'$.

\medskip
Note that ${\bf K}_M$ is not a fine moduli space, so there is no a universal
family of
pseudo-ample \pol K3 surfaces. We shall usually substitute it with a family
$f:{\cal X}\to S$ of
pseudo-ample \pol K3 surfaces (in sense of Reamrk (3.4)) such that the period
map
$S\to {\bf K}_M$ is of finite degree. A similar family $f':{\cal X}'\to S'$ of
$\check M$-polarized K3 surfaces will be called a mirror family.

\medskip
The mirror correspondence works especially nicely when $m = 1$. This is true if
and only if
$M^\perp$ contains an isotropic vector with ${\rm div}(f) = 1$. Choose $U'\cong
U$ containing $f$.
Then $M^\perp = U'\perp \check M$ and $\check M^\perp = U'\perp M$. Thus, we
can use $f\in U'$ to
define
the mirror family for both ${\bf K}_{M}$ and ${\bf K}_{\check M}$. Since
$$\check {\check M} = M,$$
we obtain
that the mirror correspondence is a duality. Note that additional assumptions
on $M$ guarantee that this duality is independent of the choice of $U'$.
For example, suppose $l(A(M) \le t-3$.
Then
$A(\check M) \cong A(M^\perp) \cong A(M)$
and by Corollary 1.13.3 of {\bf [31]} $\check M$ is determined uniquely by its
signature and the discriminant form.
Proposition 1.15.1 of loc. cit. implies that the $U$-splitting of $M^\perp$ is
unique up to an isometry of $M^\perp$.
Applying Proposition (5.1) to $M^\perp$ we obtain that $M^\perp$ admits a
unique
primitive embedding in $L$. Thus any
isometry of $M^\perp$ lifts to an isometry of $L$. This shows that the moduli
space
${\bf K}_{\check M}$ is independent of the choice of splitting $M^\perp =
U\perp \check M$.

\medskip\noindent
\ \ \ {\bf Remark (6.1).} One of the main motivation of Nikulin's paper {\bf
[31]} was
to find some conditions ensuring that
two lattices $S$ and $S'$ are K3-dual, i.e., can be realized as the lattices
$M$ and
$\check M$ from above.
For example, he proves in {\bf [31]} (Corollary 1.13.5) that two hyperbolic
lattices $S$ and $S'$ are K3-dual if ${\rm rank}\  S+{\rm rank}\  S' = 20$ and
$A(S)\cong A(S')$ as abelian groups and the values of the discriminant
quadratic forms differ
by sign. The notion of K3-dual lattices plays an important r\^ole in the
explanation of the
Arnold's Strange Duality where $M$ occurs as the lattice generated by algebraic
cycles
supported at infinity for a K3-smoothing of one of the fourteen unimodal
exceptional singularities and $\check M\perp U$ is realized
as the Milnor lattice of vanishing cycles for the same singularity. The Strange
Duality switches
the role of the lattices $M$ and $\check M$.
In {\bf
[31]} Nikulin
proves that
the Milnor lattice of a hypersurface surface singularity
contains an $1$-admissible isotropic vector whenever it is indefinite.

\medskip
For any K3 surface $X$
we can introduce the tube domain (the {\it Picard tube domain})
$${\rm Ptd}(X) = {\rm Pic}(X)_\bfr+iC(X)^+ .$$
Now let $(X,j)$ be an $M$-polarized K3 surface and $f\in M^\perp$ be an
$m$-admissible
isotropic vector. Fix a splitting
$M^\perp= U'\perp \check M$ where $U'\cong U(m)$ and $f\in U'$.
Let us consider the tube domain
$\Omega_f = V_f+iV_f^+$.
Observe that
$$V_f = ((\bfz f)^\perp_{M^\perp}/\bfz f)_\bfr \cong \check M_\bfr \subset
L_\bfr.$$
Let us choose the component $V_f^+$ such that under the above isomorphism
$$V_f^+ = V(\check M)^+ .$$
Let
$$V_f^{++} = C(\check M)^+_\bfr = \{y\in V(\check M)^+: (y,\delta) >
0\quad\hbox{for all $\delta\in \Delta(\check M)$}\},$$
$$\Omega_f^+ = V_f+iV_f^{++} = \check M_\bfr+C(\check M)^+_\bfr.$$
Let $(X,\phi)$ be an ample $\check M$-polarized surface. Then
$$\Delta(\check M)^+ = j_\phi^{-1}(\Delta(X)^+).$$
The map
$j_\phi:\check M\to {\rm Pic}(X)$ defines an open subset
$$V_f^{++}(X,\phi)= j_\phi^{-1}(C(X)^+)$$
of $V_f^{++}$, and a holomorphic embedding
$$V_f+iV_f^{++}(X,\phi) \hookrightarrow {\rm Ptd}(X).$$
Note that, if $j_\phi$ is an isomorphism we get
$V_f^{++}(X,\phi) = V_f^{++}$ and the previous embedding becomes an
isomorphism.

Let $g\in U'_\bfr$ be an isotropic vector with $(f,g) = 1$. By Remark (4.5),
it defines an isomorphism from each connected component of
$D_M$ onto the tube domain $\Omega_f$. Let $D_M^+$ be the pre-image of
$\Omega_f^+$ under this isomorphism and let ${\cal K}_M^+$ be the pre-image of
$D_M^+$ under the period map from
Theorem (4.2). For any
ample $\check M$-polarized marked K3 surface $(X,\phi)$ with bijective
$j_\phi: \check M\to {\rm Pic}(X)$, the period map defines a holomorphic
isomorphism:
$${\cal K}_{M}^+ \cong  {\rm
Ptd}(X).$$
  Note that its definition depends on the choice
of splitting $M^\perp = U'\perp \check M$, the choice of an isotropic
vector $f\in U(m)$, and the choice of marking $\phi$.

\medskip
 Recall that the period space $D_M$ admits a compactification $D_M^*$ which is
isomorphic to the
quadric in ${\bf P}(N_\bfc)$ defined by the lattice $N$. The topological
boundary of $D_M$ in $D_M^*$ is equal to the
disjoint union of locally closed analytic subsets $F$, called the boundary
components. Each boundary component is of the form
${\bf P}(I_\bfc)\cap \overline{D_M}$ for some isotropic subspace $I$ of
$N_\bfr$. Since $N$ is of signature $(2,19-t)$, we have either
${\rm dim} I = 1$ ($F$ is a point) or ${\dim} I = 2$ ($F$ is isomorphic to
upper half plane). The
stabilizer group $N(F) = \{g\in O(N_\bfr): g(F)  = F\}$ of $F$ is a maximal
parabolic subgroup of
$G = O(N_\bfr)$. Conversely, each such subgroup occurs as $N(F)$ for some
boundary component $F$. A
boundary component $F$ is called rational if
the corresponding isotropic subspace can be defined over ${\bf Q}$. It is clear
that we can identify
the set of isotropic subspaces of
$N_{\bf Q}$ with the set of primitive isotropic sublattices of $N$. In
particular we have a bijective
correspondence
$$ \{\hbox{$0$-dimensional rational boundary components of
$D_M$}\}\longleftrightarrow
I(N).$$

Let $\Gamma\subset G({\bf Q}) = O(N_{\bf Q})$ be an arithmetic subgroup of $G$
(e.g., a subgroup of
finite index in $O(N_{\bfz})$). It acts on the set ${\cal RB}(D_M)$ of rational
boundary components of $D_M$, and for each such component $F$, the stabilizer
group
$N_\Gamma(F) = \Gamma\cap N(F)$ acts discretely on $F$ with algebraic quotient
$F/N_\Gamma(F)$. Same is true for $D/\Gamma$. We have
$$\overline{D_M/\Gamma} = D/\Gamma\coprod (\bigcup_{F\in {\cal
RB}(D_M)}F)/\Gamma = D/\Gamma\coprod (\bigcup_{F\in {\cal
RB}(D_M)/\Gamma}F/N_\Gamma(F))$$
is a normal projective algebraic variety (Baily-Borel-Satake compactification).

We shall apply it to our situation when $\Gamma = \Gamma_M$. Let $f\in I(N)$
and let $F$ be the
corresponding zero-dimensional rational boundary component of $D_M$. We set
$$Z_M(f) = \{g\in N_{\Gamma_M}(F): g(f) = f\}.$$

\medskip
Now let us assume that $f\in I(N)$ is $m$-admissible and fix a splitting
$M^\perp = U'\perp \check M$ where $U'\cong U(m), f\in U'$.
Let $g\in U'$ be an isotropic vector with $(g,f) = m$.

\medskip
\proclaim\ \ \ \ Proposition (6.2).  Let $O(\check M)^* = Ker(O(\check M)\to
O(A(\check M))$. Then there is
a canonical split extension of groups
$$0\to m\check M\to Z_{M}(f) \to O(\check M)^*\to 1.$$

{\sl Proof.} We can write any $n\in M^\perp$ in the form $n = af+bg+z$, where
$a,b\in \bfz, z\in \check M$. Any $\sigma\in Z_{M}(f)$ is defined by the
formula
$$\sigma(f) = f,\ \ \sigma(g) = -{(v_\sigma,v_\sigma)\over 2m}f+g+v_\sigma, \ \
\sigma(z) =
-{(v_\sigma,\tilde \sigma(z))\over m}f+\tilde \sigma(z),$$
for some $v_\sigma,\tilde \sigma(z)\in \check M.$ It is easy to check that
$\tilde \sigma:z\to \tilde \sigma(z)$ is an element of $O(\check M)$.
Setting $A(\sigma) = (\tilde \sigma,v_{\sigma})$ we verify that
$$A(\sigma'\circ \sigma) = (\tilde \sigma'\circ\tilde \sigma, \tilde
\sigma'(v_\sigma)+v_{\sigma'}).$$
Let $G$ be the group of pairs $(s,v)\in O(\check M)\times \check M$ with the
composition law
$(s',v')\circ (s,v) = (s'\circ s,s'(v)+v').$ The homomorphism
$(s,v)\mapsto s$ makes it an extension of $O(\check M)$ with
help of $\check M$. It splits by the section $s\mapsto (s,0).$ The homomorphism
$\sigma\mapsto A(\sigma)$ is an injective homomorphism from $Z_M(f)$ to $G$.
To find its image we have to decide which pairs $(\tilde \sigma,v_\sigma)$
correspond to
isometries $\sigma\in O(M^\perp)$ which lift to isometries from $\Gamma(M)$.
By Proposition (3.3),
the condition for this is that $\sigma\in {\rm Ker}(O(M^\perp)\to
O(A(M^\perp))$.
It is easy to check that any $\sigma$ with $A(\sigma) = (1,v_\sigma)$ satisfies
this condition if and only
if ${v_\sigma\over m}\in\check M$.
Each $\sigma$ with $A(\sigma) = (\tilde \sigma,0)$
satifies this condition if and only if $\tilde \sigma\in O(\check M)^*$. Since
$G$ is the semi-direct product of
$\check M$ and $O(\check M)$, we get that the image of $Z_{M}(f)$ is the
semi-diect product of
$m\check M$ and $O(\check M)^*$. This proves the lemma.

\medskip
Let $Z_M(f)^+$ be the subgroup of
$Z_{M}(f)$ whose image in $O(\check M)^*$ consists of
elements preserving $C(\check M)^+$. The group
$Z_{M}(f)$ acts naturally on $\Omega_f = \check M_\bfr+i\check M_\bfr^+$ by the
formula
$$(\tilde \sigma,v_\sigma) (x+iy) = x+v_\sigma+i\tilde \sigma(y).$$
The subgroup $Z_M(f)^+$ preserves the tube domain
$\Omega_f^+ =
V_f+iV_f^{++} = \check M_\bfr+iC(\check M)_\bfr^+$.
It follows from the theory of compactification of homogeneous symmetric domains
that there exists
a $N_\Gamma(F)$-invariant neighborhood $\tilde U^*$ of $F$ in $D_M^*$ such that
the map
$\tilde U^*/N_\Gamma(F)\to D_M^*/\Gamma$ is an analytic isomorphism to a
neighborhood $U^*$ of the boundary point
$F/N_\Gamma(F)$ of $\overline{D_M/\Gamma}$. Restricting this isomorphism to
$\tilde U = \tilde U^*\cap \Omega_f^+$ we obtain an
isomorphism
$$\alpha:\tilde U/Z_{M}(f)^+  \to
U_F \subset U_F^*.$$
The multi-valued map
$$\alpha^{-1}:U_F \to \tilde U_F \subset \Omega_f^+ \cong {\rm Ptd}(X')$$
with the monodromy group $Z_M(f)^+$
is the {\it mirror map} MS4$'$ mentioned in the introduction.

\medskip
\ \ \ {\bf Remark (6.3).} By the Global Torelli Theorem for algebraic K3
surfaces
the group $Z_M(f)^+/\check M$ contains a subgroup of finite index isomorphic to
the
automorphism group of any surface with
Pic$(X) \cong \check M$.

\medskip

Let $(X,j)$ be a $M$-polarized K3 surface such that its isomorphism class
$[(X,j)]$ belongs to the
open subset $U_F$ from above. The pre-image of $U_F$ in $D_M$ is equal to the
disjoint union of
$\Gamma$-translates of $\tilde U_F$. So, we can choose a marking $\phi$ of
$(X,j)$ such that
the period point $P(X,\phi)$ belongs to $\tilde U_F$. Let
$\omega$ be a holomorphic 2-form on $X$, the function
$$\psi([(X,j)]) = \int_{\phi^{-1}(f)}^{}\omega$$
is a single-valued holomorphic function on $U_F$ (because the 2-cycle
$\phi^{-1}(f)$ does not depend on the marking modulo the action of the group
$Z_M(f)^+$). Thus if we normalize $\omega$ by replacing it with
$\omega' = \omega/\int_{\phi^{-1}(f)}^{}\omega$, we will be able to choose a
representative
$\ell$ of $P(X,\phi)$ with the property
$(\ell,f) = 1$.  By Corollary (4.4) to
Theorem (4.2), we obtain that this normalization allows us to identify
the Griffiths-Yukawa quadratic form on $H^1(X,\Theta_{X})$ with
the complex quadratic form $\check M_\bfc$. In particular, it gives an integral
structure
on $H^1(X,\Theta_{X})$ compatible  with the quadratic forms. This gives MS3$'$
from introduction.

\vglue .3 in
{\bf 7. Mirrors for the family of degree $2n$ polarized K3 surfaces.}
 Here we consider the mirror construction in the case $M = <2n>$. It is known
({\bf
[17]}, Theorem 1.1.2)
that $M$ admits a primitive embedding in $L$, which is unique modulo $O(L)$.
Since the lattice $U$ represents
any even integer, we may assume that
$M\subset U$ where $L = U^{\perp 3}\perp E_8^{\perp 2}$. This immediately
implies that $M^\perp
\cong U\perp U\perp
E_8\perp E_8 \perp <-2n>$. Write
$$n = \tilde n k^2,$$
where $\tilde n$ is square free. Then the group $\Gamma_M$ has exactly
$[{k+2\over 2}]$ orbits in the set of primitive isotropic vectors
in the lattice $M^\perp$ (see {\bf [37]}, Theorem 4.01). Each orbit is
represented by a vector $f$
with div$f = d, d|k$, and
$$({\bf Z} f)^\perp/{\bf Z} f \cong U\perp E_8\perp E_8 \perp <-2N>:= M_{N},$$
where $N = n/d$. So we have $[{k+2\over 2}]$ mirror families, each one is
isomorphic to
${\bf K}_{M_{n/d}}$ for
some $d|k$. Since the lattice $M_{N}$ admits a unique embedding into $L$ up to
isometry of $L$,
the number of non-isomorphic mirror moduli spaces for ${\bf K}_M$ is equal to
the number of
divisors of $k$. To study the mirror moduli spaces we may assume that
$d = 1$ by replacing $n$ with $n' = n/k.$ The corresponding isotropic vector
can be taken from a
copy of $U$.
The mirror family ${\bf K}_{\check M}$ is one-dimensional and is isomorphic to
$D_{M_n}/\Gamma_{M_n}$. We have $(M_n)^\perp = U\perp <2n>$.
So if we choose a standard basis $(f,g)$ of $U$, then we can
find a representative $\mu$ of a point from $D_{M_n}^+$ in the form
$$\mu = -nt^2f+g+te, t\in \bfc.$$
The map $\mu \to t$ defines an isomorphism from $D_{M_n}^+$ to
$\Omega_f^+$. The latter can be identified with the
upper half-plane $H =\{t=x+iy\in\bfc: y >0\}$. The group
$Z_{\Gamma(M_n)}(F)^+$ is isomorphic to
$<2n> \cong \bfz$. Let $T$ be a generator of $Z_{\Gamma(M_n)}(F)^+$
corresponding to the generator
$e$ of $<2n>$. Then $T(\check m) = \check m$ for $\check m\in \check M = M_n$,
$T(f) = f, T(g) = g-nf+e, T(e) = e-2nf$ (see the proof of Proposition (6.2)).
Then $T^{-1}(e) = e+2nf$ and
$$(\mu,f) = 1, {1\over 2n}(\mu,e) = t.$$
 From this it follows that $T$ acts on $H$ by the formula $T(t) = t-1$.
We can choose the open set $\tilde U_F$ to be equal to $\{t=x+iy: y > r\}$ for
sufficiently large $r$.
The map $t\to q = e^{2\pi it}$ defines an isomorphism
$$\tilde U_F/Z_{\Gamma(M_n)}(F)^+\cong
U_F = \{q\in \bfc^*: |q| \le 1/r\}.$$
Choose any marking
$\phi:H_2(X,\bfz)\to L$ of a $M_n$-polarized K3 surface $(X,j)\in U_F$ such
that
the period ${\bf C}\phi(\omega)$ of $(X,\phi)$ belongs to $D_{M_n}^+$ and is
equal to
${\bf C}(-nt^2f+g+te)$. Then
$$t={{1\over 2n}\int_{\phi^{-1}(e)}^{}\omega\over \int_{\phi^{-1}(f)}^{}\omega}
.$$

 This is analogous to the situation in mirror symmetry for quintic 3-folds (see
{\bf
[26]}).
\medskip
Let us now compute the global monodromy group $\Gamma_{M_n}$.
Let
$$\Gamma_0(n) = \{\pmatrix{a&b\cr
c&d\cr}\in SL(2,{\bf Z}): n|c\}/(\pm 1)\subset \Gamma = PSL(2,{\bf Z}).$$
The element
$$F = \pmatrix{0&-1/\sqrt n\cr
\sqrt n&0\cr}\in PSL(2,\bfr)$$
is of order 2 and belongs to the normalizer of $\Gamma_0(n)$ in $PSL(2,\bfr)$.
It is called the
{\it Fricke involution}. If we add it to $\Gamma_0(2)$ we obtain a subgroup of
$PSL(2,\bfr)$ denoted by
$\Gamma_0(n)^+$. It is called the {\it Fricke modular group} of level $n$.

\medskip
\proclaim\ \ \  Theorem (7.1). Let $\Gamma_{M_n}'$ be the subgroup of index 2
of
$\Gamma_{M_n}$ which stabilizes the connected component $D_{M_n}^+\cong H$ of
$D_{M_n}$. Up to a conjugation in
$PSL(2,\bfr) = Aut(H)$, we have
$$\Gamma_{M_n}' = \Gamma_0(n)^+.$$
In particular,
$${\bf K}_{M_n}\cong H/\Gamma_0(n)^+.$$

{\sl Proof.} The group $\Gamma_{M_n}$ is isomorphic
to the
group $O(U\perp <2n>)^*$. Since $-1\not\in O(U\perp <2n>)^*$, the canonical
homomorphism
$O(U\perp <2n>)^*\to PSO(1,2)$ is injective.
It is known that the groups  $PSL(2,{\bf R}) = {\rm Aut}(H)$ and
$PSO(1,2)$ are
isomorphic. For example, we can establish such isomorphism by considering a
natural
representation of $SL(2,{\bf R})$ in the space $E$ of binary forms $\alpha
x^2+2\beta\sqrt {n}xy+ \gamma y^2$ equipped with the quadratic form
$Q(\alpha,\beta,\gamma) = 2(n\beta^2-\alpha\gamma)$ (=twice the discriminant).
This allows us to view any
$g=\pmatrix{a&b\cr
c&d\cr}\in SL(2,\bfr)$
as an isometry of $(U\perp <2n>)_\bfr$ defined by the matrix
$$A(g) = \pmatrix{a^2&-2ab\sqrt n&b^2\cr
-ac/n&ad+bc&-bd/n \cr
c^2&-2\sqrt {n}cd&d^2\cr}\in SO(1,2).$$
Here the basis $(y^2, 2\sqrt n xy, x^2)$ of $E$ corresponds to the basis
$(f,e,g)$ of $U\perp <2n>$.
The kernel of the map $A:SL(2,\bfr)\to SO(1,2), g\to A(g)$ is equal to $\{\pm
1\}$. The image of the map $A$ is subgroup $SO(1,2)'$ of index 2 of $SO(1,2)$
which
preserves a connected component of $D_{M_n}$. Note that $-1$ acts as the
identity on $D_{M_n}$, so
when we
extend
$\Gamma_{M_n}'$ by adding $-1_{M_n^\perp}$ and take the intersection with
$SO(1,2)$, we obtain a
a subgroup $\Gamma_{M_n}''$ of $SO(1,2)'$ isomorphic to $\Gamma_{M_n}'$. So we
may assume now that
$\Gamma_{M_n}' = \Gamma_{M_n}''$.
Let $\Gamma'$ be the pre-image
of $\Gamma_{M_n}'$ under the map $A$. Let us describe its elements. First of
all we use that,
for any matrix $g\in \Gamma'$, the coefficients of the matrix
$A(g)$ are integers.

Write
$$a =a_1\sqrt{a_2}, b = b_1\sqrt{b_2},\ \ c = c_1\sqrt{c_2},\ \ d =
d_1\sqrt{d_2},$$
where $a_2,b_2,c_2,d_2$ are square free. We have
$$ab\sqrt n = ka_1b_1\sqrt{\tilde n a_2b_2}\in \bfz \Longrightarrow
a_2b_2\tilde n\in \bfz
\Longrightarrow a_2b_2 = s^2\tilde n$$
for some integer $s$. Since $s$ divides $a_1^2,b_1^2$ and $ab\sqrt{n}$, it must
divide the first row of the matrix $A(g)$. This
implies that $s =\pm 1$ and hence
$$a_2b_2 = \tilde n.$$
Similarly, we obtain
$$c_2d_2 = \tilde n.$$

Now, in view of above, $ac/\sqrt n\in \bfz$ gives $\sqrt{c_2/b_2},
\sqrt{a_2/d_2}\in \bfz$. This
implies that
$$c_2 = b_2, d_2 = a_2,\  k|a_1c_1,\  k|b_1d_1.$$
Let us rewrite the matrix $A(g)$ using the previous information.
$$A(g) = \pmatrix{a_1^2a_2&-2a_1a_2b_1b_2k&b_1^2b_2\cr
-a_1c_1/k&a_1d_1a_2+b_1c_1b_2&-b_1d_1/k \cr
c_1^2b_2&-2c_1d_1c_2b_2k&d_1^2a_2\cr}$$

Next we use that the discriminant group $A(U\perp <2n>)$ is generated by
the coset of ${1\over 2n}e$ modulo $U\perp <2n>$, where $e$ generates $<2n>$.
Thus elements of $O(U\perp <2n>)^*$ send $e$ to
$e+nv$ for some $v\in U$. This implies that $ad+bc\equiv \pm 1$ mod $2n$.
Assume that
$ad+bc \equiv 1$ mod $2n$. Together with $ad-bc = 1$ this implies that $n|bc$,
and hence
$k^2a_2|c_1b_1$. If a prime $p$ divides $a_2$, it must divide $c_1$ since
$p|b_1$ implies that
$p$ divides the first row of $A(g)$. On the other hand, $p|c_1$ implies that
$p$ divides the third row of $A(g)$. Thus $a_2 = 1$, hence
$$a_2 = d_2 = 1, b_2 = c_2 = \tilde n, \ \ k^2|b_1c_1$$

Let  $p$ be a prime dividing $k$. It divides $a_1c_1$ and $b_1d_1$. Assume
$p|b_1$, then $p|c_1$
since otherwise $p$ divides the whole first row of $A(g)$. Conversely, if
$p|c_1$ then $p|b_1$. Thus
$k|c_1$ and $k|b_1$, and we get
$$g = \pmatrix{a_1&kb_1'\sqrt{\tilde n}\cr
kc_1'\sqrt{\tilde n}&d_1\cr} = \pmatrix{a'&b'\sqrt{n}\cr
c'\sqrt{n}&d'\cr},\eqno (*)$$
where $a',b',c',d'\in \bfz$.
If $ad+bc\equiv -1$\ mod $2n$, we obtain similarly that
$$g = \pmatrix{a'\sqrt n&b'\cr
c'&d'\sqrt n\cr}.\eqno (**)$$
Thus we obtain that $\Gamma'$ is equal to the subgroup of $PSL(2,\bfr$ of
matrices of form (*)
and (**). Obviously matrices of form (*) form a subgroup of index $2$ in
$\Gamma'$. The whole group is generated by this subgroup and the matrix
$g_0= \pmatrix{0&-1\cr
1&0\cr}$. Now
$$\pmatrix{1/\root 4\of {n}&0\cr
0&\root 4\of {n}\cr}\cdot \pmatrix{a'&b'\sqrt{n}\cr
c'\sqrt{n}&d'\cr} \cdot \pmatrix{1/\root 4\of {n}&0\cr
0&\root 4\of {n}\cr}^{-1} = \pmatrix{a'&b'\cr
nc'&d'\cr},$$
$$\pmatrix{1/\root 4\of {n}&0\cr
0&\root 4\of {n}\cr}\cdot \pmatrix{0&-1\cr
1&0\cr}\cdot \pmatrix{1/\root 4\of {n}&0\cr
0&\root 4\of {n}\cr}^{-1} = \pmatrix{0&-1/\sqrt n\cr
\sqrt n&0\cr}.$$
This proves the theorem.

\medskip
\ \ \ {\bf Remarks (7.2)} 1. Let use the isomorphism $\Phi:H\to D_{M_n}^+, t\to
-nt^2f+g+te$. Then
$g = \pmatrix{\alpha&\beta\cr
\gamma&\delta\cr}\in SL(2,\bfr)$ acts on $H$ by the Moebius transformation
$t\to (\alpha t+\beta)/(\gamma t+\delta)$, and
$$\Phi(g(t)) =  -n(\alpha t+\beta)^2f+(\gamma t+\delta)^2g+(\alpha
t+\beta)(\gamma t+\delta)e = $$
$$=-nt^2(\alpha^2f-{\gamma^2\over n}g-{\alpha\gamma\over n}
e)+(-n\beta^2f+\delta^2g+\beta\delta e)+t(-2n\alpha\beta f+2\gamma\delta
g+(\alpha\delta+\beta\gamma)e).$$
This shows that the transformation $\Phi\circ g\circ \Phi^{-1}$ of $D_{M_n}^+$
is defined, in
the basis $(f,e,-g)$ by the matrix
$$A'(g) = \pmatrix{\alpha^2&-2n\alpha\beta&n\beta^2\cr
-\alpha\gamma/n&\alpha\delta+\gamma\beta&\beta\delta\cr
\gamma^2/n&-2\gamma\delta&\delta^2\cr}.$$
Now if $g\in \Gamma_0(n)^+$ we observe that $A'(g)\in \Gamma_{M_n}'$. This
shows that
$$\Phi\circ \Gamma_0(n)^+\circ\Phi^{-1} = \Gamma_{M_n}'.$$

\medskip
2. It is known that the orthogonal group of the discriminant group of the
lattice
$U\perp <2n>$ is isomorphic to the group $(\bfz/2\bfz)^s$, where $s$ is the
number of
distinct prime divisors of $n$ ({\bf [37]}, Lemma 3.6.1). If $n= k^2\tilde n$
as above with
$(k,\tilde n) = 1$, this group is isomorphic to $\bar\Gamma_0(n)/\Gamma_0(n)$
where
$\bar\Gamma_0(n)$ is the abelian normalizer of $\Gamma_0(n)$ in $SL(2,\bfr)$
(see {\bf [20]}, Theorem 3).  Using Nikulin's results from {\bf [31]}, one can
show that the homomorphism $O(U\perp<2n>)\to
O(A(U\perp<2n>))$ is surjective. Since
$\Phi\circ \Gamma_0(n)\circ\Phi^{-1}$ is equal to $\Gamma_{M_n}\cap
SO(U\perp<2n>)$, this easily implies that
$$\Phi\circ \bar\Gamma_0(n)\circ\Phi^{-1} = O(U\perp<2n>).$$
The group $O(A(U\perp<2n>))$ acts on ${\bf K}_{M_n}$ with kernel isomorphic to
$\{\pm 1\}$. The
quotient is the moduli space of K3 surfaces admiting a pseudo-ample
$M_n$-polarization.

\medskip
Let us now find the subset ${\bf K}_{M_n}^a\subset H/\Gamma_0(n)^+$ of
isomorphism classes
of ample $M_n$-polarized K3 surfaces.

\medskip
\proclaim\ \ \  Theorem (7.3). Let $S\subset H/\Gamma_0(n)^+$ be the set of
orbits of the points
${c\over b}+{i\over b\sqrt n}$, where $c\in \bfz$ and $b|cn^2+1$. Then
$${\bf K}_{M_n}^a = H/\Gamma_0(n)^+\setminus S.$$
Moreover,
$$\#S = \cases{1&if $n\le  4$,\cr
2h(-4n)&if $n\equiv 7$ mod $8$,\cr
4h(-4n)/3&if $n\equiv 3$ mod $8, n\ge 4$,\cr
h(-4n)&otherwise.\cr}$$
Here $h(k)$ denotes the number of classes of primitive binary quadratic forms
of discriminant
$k$.

{\sl Proof.} Recall from section 3 that ${\bf K}_{M_n}\setminus{\bf K}_{M_n}^a$
is
equal to the set of $\Gamma_{M_n}$-orbits in $D_{M_n}$ of hyperplanes
$H_v=\{z\in D_{M_n}:(z,v) = 0\}, v\in U\perp <2n>, (v,v) = -2.$  Let us use the
isomorphism
$\Phi:H\to D_{M_n}^+, t\to -nt^2f+g+te$. Let $v= af+bg+ce$ with $(v,v) =
2ab+2nc^2=-2$. Then
$(-t^2f+g+te,v) = -nbt^2+a+2nct = 0$ implies
$$t = {c\over b}\pm\sqrt {{c^2\over b^2}+{a\over nb}} ={c\over b}\pm\sqrt
{{nc^2+ab\over nb^2}} =
{c\over b}+{i\over b\sqrt n}.\eqno (*)$$
This proves our first assertion. Let $t\in H$ such that $\Phi(t)\in H(v)$ for
some hyperplane $H(v)$ as above.
Since $h(v)$ is fixed by an automorphism of order 2 corresponding to the
reflection isometry $x\to x+(x,v)v$, we see that
$t$ is fixed by some involution $g\in \Gamma_0(n)^+$. Let $g$ be represented by
a matrix $\tilde g = \pmatrix{a&b\cr
c&d\cr}.$ We have either $\tilde g^2 = 1$ or $\tilde g^2 = -1$. Since the
characteristic polynomial of
$\tilde g$ is equal to $X^2-(a+d)X+1$, we see that only the second case occurs,
and $a+d = 0$. If
$g\in \Gamma_0(n)$, then
$\tilde g = \pmatrix{a&b\cr
nc&-a\cr}$ where $a,b,c,d\in \bfz$. The fixed points $t$ of $g$ can be
computed, and we find that
$t = {a\over cn}+{i\over cn}$. This differs from points (*) unless $n=1$. If
$n=1$,
$\Gamma_0(1)^+ =\Gamma$, and there is only
one orbit of
such points. If $n > 1$, and $g$ is an involution  from the coset
$F\cdot\Gamma_0(n)$ of the Fricke
involution $F$,
we find that its fixed points look like (*). Consider the double cover
$p:X_0(n)\to X_0(n)^+ = X_0(n)/(T)$, where $X_0(n)$ (resp. $X_0(n)^+$) is a
nonsingular projective model of the quotient
$H/\Gamma_0(n)$ (resp. $H/\Gamma_0(n)^+$). We have a bijective correspondence
between the  ramification
points of this cover and $\Gamma_0(n)$-orbits in $\bar H = H\cup \{\infty\}\cup
{\bf Q}$ whose
stabilizer
belongs to the coset $F\cdot\Gamma_0(n)$. When $n\ge 5$ one checks that
$F\cdot\Gamma_0(n)$ does not have
parabolic elements (i.e.. elements which fixes
$t\in {\bf Q}\cup\{\infty\}$) and  elements of finite order greater than 2.
This shows that
$\#S$ is equal to the number of ramification points of the double cover $p$.
This number was computed
by R. Fricke in {\bf [12]}, and it is equal to the number which we gave in the
statement of the
theorem. Now, it is known that the modular curve $X_0(n)$ is of genus 0 when $n
= 2,3,4$. Thus there
are only 2 ramification points. One of them is an orbit with stabilizer of
order 2 contained in $\Gamma_0(2)$.
Another one is an orbit of with stabilizer of order 2 whose generator belongs
to $F\cdot\Gamma_0(n).$ This
proves the assertion.

\medskip
Assume now that the curve $X_0(n)^+ = \overline {H/\Gamma_0(n)^+}$ is rational.
All such $n$
can be listed (see {\bf [18]}) (as was observed by A. Ogg, the primes from this
list are just
those which divide
the order of the Fisher-Griess Monster group).
Let
$$C(n) = \overline {H/\Gamma_0(n)^+}\setminus (H/\Gamma_0(n)^+)$$
be the set of cusp points.

\medskip
\proclaim \ \ \ Proposition (7.3). Let $\phi$ be the Euler function. Then
$$\#C(n) = \cases{{1\over 2}\sum\limits_{d|n,d> 0}\phi((d,n/d))&if $n\ne 4$\cr
2&if $n = 4$.\cr}$$

{\sl Proof.} The number of cusps for the modular curve $X_0(n)$ is equal to
$\sum\limits_{d|n,d> 0}\phi((d,n/d))$ (see {\bf [38]}, Proposition 1.4.1). It
is known that the Fricke involution acts on this set without fixed points if
$n\ne 4$ (see {\bf [18]}) and has one fixed point if $n=4$. From this the
result follows.

\medskip
\proclaim \ \ \ Corollary (7.4). Let $n = p$ be a prime number, $M= <2p>$.
Assume that $X_0(p)^+$ is rational.
Then
$${\bf K}_{\check M}\cong {\bf A}^1.$$

\medskip
\proclaim \ \ \ Theorem (7.5). Assume $X_0(n)^+$ is rational. There exists a
unique holomorphic function
(called the {\it Hauptmodule})
$$j_n:H\to \bfc$$
satisfying the following conditions:
\item {(i)} $j_n$ is invariant with respect to $\Gamma_0(n)^+$;
\item{(ii)} $j_n$ has a Fourier expansion
$$j_n(t) = q^{-1}+\sum_{m=1}^\infty c_mq^m, \ \ \ q = e^{2\pi i t};$$
\item{(iii)} the coefficients of the Fourier expansion are all integers;
\item{(iv)} considered as a meromorphic function on $X_0(n)^+$, the function
$j_n$ has a simple
point at the cusp
$\Gamma_0(n)^+\cdot\infty$ and generates the field of meromorphic functions on
$X_0(n)^+$.

{\sl Proof.} See {\bf [18]}.

\medskip
Let us restrict the meromorphic function $j_n^{-1}$ to a neihborhood $\tilde
U_F = \{t=x+iy\in H:y >r\}$
for sufficiently large $r$ chosen so that $j_n^{-1}$ is holomorphic on $\tilde
U_F$. Then
the properties of $j_n$ assure that
$j_n^{-1}$ defines an isomorphism from $\tilde U_F/(\Gamma_0(n)^+)_\infty$ to a
neihborhood of the
cusp $\Gamma_0(n)^+\cdot \infty$. Comparing it with the discussion in the
beginning of the section,
we find that the mirror map at the cusp can be given by the inverse of the
Hauptmodule function
$j_n^{-1}$.
 This should be compared to ${\bf [21]}$.

\medskip
It is well-known that $H/\Gamma_0(n)$ is a coarse moduli space for the
isomorphism classes of pairs $(E,A)$, where $E$ is an elliptic curve, and
$A$ is a cyclic subgroup of order n of of $E$. The Fricke involution acts on
$H/\Gamma_0(n)$ by sending the pair $(E,A)$ to the pair
$(E/A, E_n/A)$. Let us give an explicit geometric relationship between the
isomorphism class of a $M_n$-polarized K3 surface respesented by a point
$z\in H/\Gamma_0(n)^+$ and the isomorphism class of the pair of isogeneous
elliptic curves
$(E,E'= E/A)$ represented by the same
point $z$.  This can be used to explain the observation of B. Lian and S. Yau
that the periods of certain one-dimensional
families of K3 surfaces can be expressed as the products of periods of some
family of
elliptic curves
(see {\bf [22]}). I am grateful to Dan Burns who suggested that our K3 surfaces
should be related to Kummer surfaces
Kum$(E\times E')$.

\medskip
\proclaim\ \ \  Theorem (7.6). Let $M = <2n>$ and $X$ be a $M_n$-polarized K3
surface with period $t\in H$. Let
$E_t = \bfc/\bfz+t\bfz$ and $E_t' = \bfc/\bfz+(-1/nt)\bfz$ be the corresponding
pair of
isogeneous elliptic curves. Then there exists a canonical involution $\tau$ on
$X$ such that $X/(\tau)$ is
birationally isomorphic to the Kummer surface $E_t\times E_t'/(\pm 1)$.

{\sl Proof.} The fact that there exists an involution $\tau$ on $X$ such that
$X/(\tau)$ is
birationally isomorphic to some Kummer surface $A/\{\pm 1\}$ follows from the
property that $rk M_n = 19$ {\bf [25]}. As is explained in loc.cit. and in {\bf
[30]} such an involution corresponds to a
primitive embedding $i:E_8(2)\to M_n$ (the image is the sublattice of
$\tau$-antiinvariant divisor classes). Here $E_8(2)$ denote the lattice
obtained from the lattice $E_8$ by
multiplying its quadratic form by 2. We define this embedding to be the
canonical one:
$i:E_8(2)\to E_8\perp E_8\perp U\perp <-2n>, x\to (x,x,0,0).$ Then it is shown
that
$X/(\tau)\cong {\rm Kum}(A) = A/\{\pm 1\}$, where
$A$ is an abelian surface. Let $Y$ be a minimal nonsingular model of Kum$(A)$.
The rational map
$\pi:X\to Y$ induces an embedding of lattices of trancendental cycles
$\pi^*:T_Y(2) \to T_X$. It is also known {\bf [32]} that $\pi^*(T_Y(2)) = 2S$
where
$S\subset T_X\otimes {\bf Q}$ with
$S/T_X \cong (\bfz/2\bfz)^\alpha \subset A(T_X)$. If $X$ satisfies $Pic(X) =
M_n$, then
$T_X = U\perp <2n>$, and it is easy to see that $T_Y = U(2)\perp <4n> =
T_X(2)$. Also, it is known that
$T_A(2) \cong T_Y$ (see [4], Chapter VIII, \S5).  Let $p:A\to Y$ be the
rational map of degree 2
defined by the canonical map $A\to {\rm Kum}(A)$. It follows from loc.cit. that
the homomorphism
$p^*:H^2(Y,\bfc)\to H^2(A,\bfc)$ preserves the Hodge structures, i.e.,
$p^*(H^{2,0}(Y)) = H^{2,0}(A)$. The same property is true
for $\pi_\bfc^*: (T_Y)_\bfc\to (T_X)_\bfc$. Thus the isomorphism
$p^*\circ(\pi^*)^{-1}(T_X)_\bfc\to (T_A)_\bfc$ preserves
the Hodge structures. So let us compute the period of the abelian surface $A$,
knowing that the
period of $X$ is equal to $\mu = -nt^2f+g+te\in D^+_{M_n}$. Recall that for any
complex torus
$T$, we have an isomorphism of lattices
$H^2(T,\bfz)\cong U\perp U\perp U = U^{\perp 3}$. Fix a primitive embedding
$i:U\perp <2n>\hookrightarrow U^{\perp 3}$. To be more precise, let $\Lambda =
\bfz e_1+\bfz e_2+\bfz e_3+\bfz e_4$ with a
fixed isomorphism
$d:\bigwedge^4 \Lambda \to \bfz$ such that $d(e_1\wedge e_2\wedge e_3\wedge
e_4) = 1$. Then $\bigwedge^2 \Lambda$ has a structure of a lattice with respect
to the bilinear form
$(\alpha,\beta) = d(\alpha\wedge\beta)$. Consider the following basis of
$\bigwedge^2\Lambda$:
$$f_1 = e_1\wedge e_2,\  g_1 = e_3\wedge e_4,\  f_2 = e_1\wedge e_3,\  g_2 =
e_4\wedge e_2,\  f_3 = e_1\wedge e_4,
\ g_3 = e_2\wedge e_3.$$
Then $\bfz f_i+\bfz g_i \cong U$, and $\bfz f_i+\bfz g_i$ is orthogonal to
$\bfz f_i+\bfz g_i$ for $i\ne j$. Our embedding $i:U\perp <2n>\to U^3$ can be
chosen as follows:
$i(f) = f_1, i(g) = g_1, i(e) = f_2+ng_2$. For simplicity of notation we denote
by
$i:(U\perp <2n>)_\bfc\to
(U^{\perp 3})_\bfc$ the extension of the embedding $i$ to the injective map of
the complexified
spaces. We have
$$i(\mu) = -nt^2f_1+g_1+t(f_2+ng_2) = -nt^2e_1\wedge e_2+e_3\wedge
e_4+t(e_1\wedge e_3+ne_4\wedge e_2).$$
We immediately verify that
$$i(\mu) = (-te_1+e_4)\wedge(nte_2-e_3).$$
Using {\bf [13]}, Expos\'e VIII, we can interpret it as follows. Let
$$E_t = \bfc/\bfz+t\bfz,\ \  E_t' =
\bfc/\bfz+(-1/nt)\bfz$$
be
the  pair of isogeneous elliptic curves. Then under a certain marking
$\phi:H^2(E_t\times E_t',\bfz)\to U^{\perp 3}$, $\phi(H^{2,0}(E_t\times E_t'))
= i(\mu).$ Now the
assertion follows from the Global
Torelli Theorem.

\medskip

Let us exhibit explicitly some mirror families of the family of polarized K3
surfaces of degree $2n$.
We shall use the notation
$A_n,D_n,E_n$ to denote the negative definite even lattice defined
by the negative of the Cartan matrix of the root system of a simple Lie algebra
of type
$A_n,D_n,E_n$, respectively. We shall use the following well-known description
of the Picard lattice of an
elliptic surface $f:X\to S$ with a section:

\medskip\noindent
\proclaim\ \ \  Lemma (7.7)(Shioda-Tate). Let $Pic(X)'$ be the subgroup of
$Pic(X)$ generated by
irreducible components of
fibres and by a section. Then the quotient group $Pic(X)/Pic(X)'$ is isomorphic
to the
Mordell-Weil group $MW(X/S)$ of
sections of the fibration.

{\sl Proof.} See {\bf [7]}, Proposition 5.3.4.

\medskip
This lemma is applied as follows.  We exhibit an elliptic fibration such that
$Pic(X)'$ is a
subgroup of finite index in $Pic(X)$. The lemma implies that the Mordell-Weil
group is finite. Then we show
that it is in fact trivial. Now it is easy to find the structure of the lattice
$Pic(X)'$. Its
sublattice generated by a section and a fibre is isomorphic to the lattice
$U$. Its orthogonal complement is isomorphic to the sum of lattices of type
$A_n, D_n,E_n$, each is spanned by the
irreducible components of a fibre which do not intersect the chosen section.

\medskip\noindent
\ \ \ {\bf Example (7.8) ($n = 1$).}
We have
$$M_n \cong U\perp E_8\perp E_8\perp\langle -2\rangle.$$
Using the previous remark, it suggests to look for a K3 surface with an
elliptic fibration
$f:X\to {\bf P}^1$ with a section and two reducible fibres of types $\tilde
E_8$ (or $II^*$ in Kodaira's notation),
and one reducible fibre of type $\tilde A_1$ (Kodaira's $I_2$ or $III$).
Since the group $F^\natural$ of non-singular points of a
fibre $F$ of type $\tilde E_8$ is isomorphic the additive group $\bfc$, and the
restriction homomorphism MW$(X/{\bf P}^1)\to F^\natural$ is known to be
injective on the
torsion subgroup ({\bf [7]}, Proposition 5.3.4), we obtain that Tors(MW$(X/{\bf
P}^1))$ is trivial. Hence, if $Pic(X)$ is known to
be of rank $19$, it must be isomorphic to $M_n$.

To
construct such a surface $X$, we
take a nonsingular plane cubic $C$ and the tangent line $L$ at its inflection
point. The pencil of
plane cubics spanned by $C$ and $3L$ defines a rational map
${\bf P}^2\to {\bf P}^1$. In appropriate coordinate system we can represent the
pencil in the form:
$$\lambda(Y^2Z+X^3+aXZ^2)+\mu Z^3 = 0.$$
 After resolving its nine fundamental points (infinitely near to the point
$(0,1,0)$)
we arrive at a
rational elliptic surface $f:V(a)\to {\bf P}^1$. It contains
a degenerate fibre of type
$\tilde E_8$ corresponding to $(\lambda,\mu) = (0,1)$. The irreducible singular
fibres correspond
to $(\lambda,\mu) = (1,b)$ where $4a^3+27b^2 = 0$. If $a\ne 0$ we have two
irreducible
singular fibres with ordinary double points. If $a = 0$ we have one irreducible
singular fibre with a cusp singularity.
Let $F_a(b)$ denote the fibre of $f$ corresponding to $(\lambda,\mu) = (1,b)$.
Let $F_1 =
F_a(b)$, where $4a^3+27b^2 = 0,$ be an irreducible singular fibre, and let
$F_2 = F_a(b+1)$.
Consider the
double cover $X(a,b)'$ of $V(a)$ branched along the union $F_1\cup F_2$. After
resolving its singularities we
obtain a K3 surface $X(a,b)$
with two reducible fibres of type $\tilde E_8$. It has additional reducible
fibres:
one fibre of type $\tilde A_1$ if $F_1$ has a node, $F_2$ is nonsingular
($a\ne 0,b\ne -{1\over 2}$); two fibres of
type  $\tilde A_1$ if $F_1, F_2$ have nodes ($b = -{1\over 2}$); one fibre
of type
$\tilde A_2^*$ (Kodaira's IV) if $F_1$ has a cusp ($a = 0$).

We have a one-parameter family ${\cal X}'\to C$ of singular surfaces $X(a,b)'$
parametrized by the
affine curve $C:4A^3+27B^2 = 0$. The map ${\cal X}'\to C$ is equivariant with
respect to the natural
action of the group $\mu_3$ of
third roots of unity. Its generator $\rho = e^{2\pi i/3}$ acts on
$C$ by $(a,b)\to (\rho a,b)$ and on ${\cal X}'$
via its action
on ${\bf P}^2$ by the formula $X\to \rho X, Y\to Y, Z\to Z$. After dividing
${\cal X}'$ by this
action, we obtain a family
$\pi':{\cal Y}' = {\cal X}'/\mu_3\to {\bf A}^1 = C/\mu_3$. Let $0\in {\bf A}^1$
be the orbit of $(0,0)\in C$ and $1\in {\bf A}^1$ be the orbit of
$(a,-{1\over 2})$. For any $t\in {\bf A}^1\setminus\{0,1\}$, the fibre  ${\cal
Y}'_t= \pi'^{-1}(t)$
has one ordinary double
point. The fibre
${\cal Y}'_0\cong X(0,0)/\mu_3$, it is a rational singular surface. The fibre
${\cal Y}_1 \cong X(a,-{1\over 2})$. It has two ordinary double points. Let
$$\pi:{\cal Y}\to {\bf A}^1\setminus\{0,1\}$$
be the composition of $\pi'$ and the blowing up ${\cal Y}\to \pi'^{-1}({\bf
A}^1\subset\{0,1\}$ of the
locus of singular points of the fibres ${\cal Y}'_t, t\ne 0,1$. We have
constructed a family of
pseudo-ample $M_1$-polarized K3 surfaces.
The period map for the family ${\cal Y}\to {\bf A}^1\setminus \{0,1\}$ defines
a regular map:
$$p: {\bf A}^1\setminus\{0,1\}\to {\bf K}_{M_1}\cong H/\Gamma\cong {\bf A}^1.$$

One can show that the period mapping $p$
can be extended to an isomorphism ${\bf A}^1 \to H/\Gamma$ which sends the
point $0$ to the
orbit of $e^{2\pi i/3}$ and $1$ to the orbit of $i$. The first point is a
period of the surface $X(0,0)$, the
second point is a period of the surface $X(a,-{1\over 2})$. The latter surface
is non-ample
$M$-polarized surface. The monodromy group of our family is generated by the
local monodromies at
$0,1$ and $\infty$. They are isomorphic to the subgroups of $\Gamma$ which
stabilize
$e^{2\pi i/3},i$ and $\infty$, respectively. Thus the global monodromy of our
family is isomorphic to
$\Gamma$.

Since not all fibres are isomorphic,
the period map is not constant.
Hence there exists a dense subset $U$ of ${\bf A}^1\setminus \{0,1\}$ such that
the Picard number of
${\cal Y}_t, t\in U,$ is equal to 19, and hence Pic${\cal Y}_t \cong M_1$.

Observe that $X(0,0)\cong X_0$ has the Picard lattice of rank 20 isomorphic to
$U\perp E_8\perp E_8\perp A_2$, and $X(a,-{1\over 2})$ has Picard lattice of
rank 20
isomorphic to $U\perp E_8\perp E_8\perp <-2>\perp <-2>.$

\medskip\noindent
\ \ \ {\bf Example (7.9) ($n = 2$).}
 So we want to describe a mirror family for quartic surfaces.
One can show, for example, using the uniqueness results from {\bf [31]}, that
$$M_n = U\perp E_8\perp E_8\perp <-4> \cong  U\perp E_8\perp D_9.$$
Similar to the previous example, we should construct a one-dimensional family
of elliptic K3 surfaces $X$ with a section, one
reducible fibre of type $\tilde D_9$ (Kodaira's $I_5^*$) and one reducible
fibre of type
$\tilde E_8$. To construct the family ${\cal F}$ of such elliptic surfaces we
use the same idea as in
the previous example. Consider the pencil of cubic curves:
$$F(\lambda,\mu) = \lambda X^3+\mu Z(Y^2-XZ+aX^2) = 0.$$
Let $V\to {\bf P}^1$ be the associated rational elliptic surface.
It has a degenerate fibre of type $\tilde E_7$ (Kodaira's $III^*$) and a
reducible fibre $F_1$ of
type $\tilde A_1(a\ne 0)$ or $\tilde A_1^* (a = 0)$.
 Let $X$ be the double
cover of $V$ branched along the union of $F_1$ and another irreducible fibre
$F_2$.
$X$ is an elliptic surface with two reducible fibres of type
$\tilde E_7$ and a reducible fibre of type
$\tilde A_3$ (or $\tilde D_4$). If $F_2$ is singular (this happens when $\mu
a^2= 4\lambda$) it has an additional  fibre of type $\tilde A_1$.
The elliptic fibration has also two sections. We claim that $X$ has another
elliptic fibration with two reducible fibres of type
$\tilde E_8$ and $\tilde D_9$. To see this we assume for simplicity that $F$ is
of type $\tilde A_3$.
Let
$$D = 2R_0+R_1+2R_2+3R_3+4R_4+3R_5+2R_6+R_7$$
$$D' = 2R_0'+R_1'+2R_2'+3R_3'+4R_4'+3R_5'+2R_6'+R_7'$$
be the reducible fibres of type $\tilde E_7$, let $F = E_0+E_1+E_2+E_3$ be the
other
reducible fibre. Without loss of generality we may assume that the two sections
$S_1$ and $S_2$
intersect the fibre $D$ at $R_1$ and $R_7$, respectively, and the fibre $D'$ at
$R_1'$ and $R_7'$, respectively.
Also $S_1$ intersects $F$ at $E_0$ and $S_2$ intersects $F$ at $E_2$. Now
consider the following
disjoint curves with self-intersection $0$:
$$D_1 = 3R_0+2R_2+3R_3+6R_4+5R_5+4R_6+3R_7+2S_2+R_7',$$
$$D_2 = R_5'+R_0'+2R_4'+2R_3'+2R_2'+2R_1'+2S_1+2E_0+E_1+E_3.$$
By Hodge's Index Theorem, the divisors $D_1$ and $D_2$ are linearly equivalent.
They span a pencil which
defines an elliptic fibration with fibre $D_1$ of type $\tilde E_8$ and fibre
$D_2$ of type
$\tilde D_9$.

Let $X(t;a)$ be the elliptic surface obtained by the above construction when we
take
$F_2 = F(\lambda,\mu)$ with $t =\lambda/\mu\ne 0,\infty$. The linear
substitution $X\to cX, Y\to Y, Z\to c^{-1}Z$
extends to an isomorphism
$X(t;a)\cong X(c^4t,c^2a)$. Let the group $\bfc^*$ act on $\bfc^*\times
\bfc^*\setminus \{(t,a):a^2 =t\}$ by the formula
$(t,a)\to (c^4t,c^2a)$. The orbit space is isomorphic to
$\bfp^1\setminus\{0,1,\infty\}$.
 When we vary $(t,a)\in \bfc^*\times \bfc$,
we obtain a
family ${\cal Y}\to \bfc^*\times \bfc$ of
$M_2$-polarized K3 surfaces with $X(t;a)\cong {\cal Y}_{(t,a)}$ for
$a^2 \ne 4t,0$. When $a^2 = 4t,0$, the fibre ${\cal Y}_{(t,a)}$ is singular but
birationally isomorphic to
$X(t;a)$. The surface $X(t;0)$ has a fibre of type $\tilde D_4$ and its Picard
number equals 20.
The surface $X(a^2/4,a), a\ne 0,$ has a reducible fibre of type $\tilde A_1$
and its Picard lattice is isomorphic to
$U\perp E_8\perp D_9\perp <-2>$. Let $f:\bfc^*\times\bfc \to \bfc$ given by the
formula $(t,a)\to a^2/4t$.
As in the previous example, we can descend the family ${\cal Y}\to \bfc^*\times
\bfc$ to a family
$\bar {\cal Y}\to \bfc$ of pseudo-ample
$M_2$-polarized K3 surfaces with singular fibres over $0= f(t,0)$ and
$1=f(a^2/4,a)$. The period map extends to an
isomorphism ${\bf A}^1 \to {\bf K}_{M_2}$ which sends $0$ to the isomorphism
class of the
surface $X(0)$ and sends $4$ to the isomorphism class of the surface $X(2)$.
The latter surface is a pseudo-ample but not ample $M$-polarized K3 surface.

\medskip\noindent
\ \ \ {\bf Example (7.10) ($n = 3$)}  We skip the details. We have
$$M_3 \cong U\perp E_8\perp E_8\perp <-6>.$$
We consider a
rational elliptic surface $V$ with a section, one reducible fibre $F_1$ of type
$\tilde E_6$ and one reducible fibre $F_2$ of type $\tilde A_2$.
To construct such a surface we take a plane nonsingular cubic $C$ and three
inflection
points on it lying on a line  (this means that they add up to 0 in the group
law on the cubic with an inflection point
taken as the origin). Then we take the pencil of cubic curves spanned by $C$
and the union of the tangent lines at the three inflection points.
After resolving the base points of this pencil we arrive at the surface $V$.
The surface $X$ is obtained as a
minimal nonsingular model of the double cover of $V$ branched over $F_2$ and a
nonsingular fibre.
The surface $X$ is an elliptic K3 surfaces with a section, two reducible fibres
of type
$\tilde E_6$ (Kodaira's $IV$) and one reducible fibre of type $\tilde A_5$
(Kodaira's $I_5$).
Its Mordell-Weyl group is $\bfz/3$ and the sublattice of Pic$(X)$ spanned by
the section and components of the reducible fibres is isomorphic to
$E_6\perp E_6\perp A_5$. By Lemma (7.4), we get that
Pic$(X)$ is a hyperbolic lattice of discriminant $6$. One can find another
pencil on this surface
with three reducible fibres of
type
$\tilde E_8, \tilde E_7$, and $\tilde A_2$. Since its discriminant equals 6, it
must coincide with $\pic X$.
On the other hand, it has the same discriminant group as the lattice $M_3$.
By Nikulin's uniqueness results, we conclude that $\pic X \cong M_3$.

\medskip
{\bf Remark (7.11).}
The following remark may be appropriate. As we have already noticed in the
previous section, the
rational
one-dimensional boundary components of $D_M$ correspond to rank 2 primitive
isotropic
sublattices $S$ of $N = M^\perp$. Each component contains in its closure the
0-dimensional boundary component
defined by an isotropic vector $f\in S$. If $N = U\perp M_n$ and $f\in
U$,
then $S$ is determined by a primitive isotropic vector in $M_n$. Now if
identify $M_n$ with ${\rm Pic}(X')$ for some ample $M_n$-polarized K3 surface
$X'$ from
the mirror moduli space,
we find a bijection between one-dimensional rational boundary components
of $D_M/\Gamma_M$ containing the given $0$-dimensional boundary component and
isomorphism classes of elliptic fibrations on  $X'$. When $M = <2>$ or $<4>$
the list of
$\Gamma_M$-orbits of two-dimensional isotropic sublattices of $N$ containing a
given primitive isotropic vector is given in {\bf [37]}. We find that
$S_N^\perp/S$ can be isomorphic to one of the following lattices:
$$A_1\perp E_8\perp E_8, A_1\perp D_{16}, E_7\perp D_{10}, A_{17} \ \ (M =
<2>).$$
$$E_8^2\perp<-4>, D_{16}\perp <-4>, E_8\perp D_9, E_7^2\perp A_3, D_{17},
D_{12}\perp D_5,$$
$$D_8^2\perp <-4>, A_{15}\perp A_1^2, E_6\perp A_{11}\ \ ( M = <4>).$$
In
our interpretation we obtain that the mirror surfaces contain elliptic
fibrations with reducible fibres of type:
$$\tilde A_1, \tilde E_8, \tilde E_8;\ \  \tilde A_1, \tilde D_{16};\ \ \
\tilde E_7,\tilde D_{10};\ \ \ \tilde A_{17} \ \ (M = <2>).$$
$$\tilde E_8,\tilde E_8;\ \ \ \tilde D_{16};\ \ \ \tilde E_8,\tilde D_9;\ \ \
\tilde E_7,\tilde E_7,\tilde A_3;
\ \ \tilde D_{17};$$
$$\ \ \ \tilde D_{12},\tilde
D_5;\ \ \
\tilde D_8,\tilde D_8;\ \ \ \tilde A_{15},\tilde A_1,\tilde A_1;\ \ \ \tilde
E_6,\tilde A_{11}\ \ ( M = <4>).$$
We have seen already a pencil of type $\tilde E_8+\tilde E_8+\tilde A_1$ on
surfaces from the mirror family of ${\bf K}_{<2>}$ and the pencils of type
$\tilde A_1+\tilde A_1+\tilde A_{15}$ and $\tilde E_8+\tilde D_{9}$ on surfaces
from the mirror family of ${\bf
K}_{<4>}$.

Similar computation is known for the case $M = <6>$ (see {\bf [39]}). We have
the following types of elliptic fibrations on surfaces from
${\bf K}_{M_3}$:
$$\tilde E_8,\tilde E_8;\ \ \ \tilde D_{16};\ \ \ \tilde E_8,\tilde E_7,\tilde
A_2;
\ \ \ \tilde D_{14},\tilde A_2,\tilde A_1;\ \ \ \tilde D_{10},\tilde D_{6};\ \
\
\ \ \tilde D_{8},\tilde E_{7},\tilde A_{1};$$
$$\tilde A_{15},\ \ \ \tilde E_{6},\tilde
E_6,\tilde A_5;\ \ \
\tilde A_{11},\tilde D_5,\tilde A_1;\ \ \ \tilde A_{9},\tilde D_7.$$
In example (7.10) we have seen an elliptic fibration of type $\tilde E_6,\tilde
E_6,\tilde A_5$.

\vglue .3 in
{\bf 8. Toric hypersurfaces.}
Recall the following mirror construction of Batyrev {\bf [5]} which generalizes
the original construction of Green-Plesser.
 Let $\Delta \subset \bfr^n$ be a convex
$n$-dimensional lattice polytope given by inequalities:
$$\sum_{j=1}^na_{ij}x_j \le 1 , i = 1,\ldots,k,$$
where $a_{ij}\in\bfz$ (a reflexive polytope).  Let $\Delta^*$ be
the polytope equal to the convex hull of the vectors $l_i =
(a_{i1},\ldots,a_{in})$. It is also a reflexive polytope. Let ${\bf P}_\Delta$
(resp. ${\bf P}_{\Delta^*}$) be the corresponding toric
variety, and ${\cal F}(\Delta)$ (resp. ${\cal F}(\Delta^*)$) be the family of
hypersurfaces in ${\bf P}_\Delta$ (resp. ${\bf P}_{\Delta^*})$ defined by
$\Delta$-nondegenerate (resp. $\Delta^*$-nondegenerate) Laurent polynomials.
For $n\le 4$ there exists a map $f: \tilde\bfp_\Delta\to \bfp_\Delta$ such that
the proper transform of a
general member of the family ${\cal F}(\Delta)$ is a Calabi-Yau manifold.
A similar construction with $\Delta^*$ defines another family of Calabi-Yau
manifolds.
In the case $n =4$ the two families of
Calabi-Yau 3-folds satisfy the first attribute of mirror symmetry: the
dimension
of the local moduli space for a member of the first family is equal to the
Picard number of
a member of the second family {\bf [5]}.

 Consider the special case of Batyrev's construction when $\Delta$ is a
3-dimensional
simplex:
$$\Delta(w) = \{(t_0,t_1,t_2,t_3)\in \bfr^4:\sum_{i=0}^3w_it_i = 0, t_i\ge -1,
i= 0,\ldots,3\},$$
where $w = (w_0,w_1,w_2,w_3)$ is a collection of four positive integers with
greatest common
divisor equal to 1
and such that $d = w_0+w_1+w_2+w_3$ is divisible by each $w_i$. Here we
identify $\bfr^3$ with the hyperplane $\sum_{i=0}^3w_it_i = 0$. The toric space
$\bfp_\Delta$
is the weighted projective space $\bfp(w) = \bfp(w_0,w_1,w_2,w_3)$.
The family ${\cal F}(\Delta)$ is the family of quasi-smooth hypersurfaces of
degree $d$ in
$\bfp(w)$. One of its representatives is the surface
$$
 x_0^{d_0}+x_1^{d_1}+x_2^{d_2}+x_3^{d_3} = 0,
$$
where $d_i = d/w_i, i = 0,\ldots,3.$
Let $\Pi$ be the finite abelian group of order $d_0d_1d_2d_3/d^2$ equal to the
kernel of the homomorphism
$$(\mu_{d_0}\times \mu_{d_1}\times \mu_{d_2}\times \mu_{d_3})/\mu_d \to \mu_d,\
\  g_0^{a_0}\cdots g_3^{a_n} \to g^{w_0a_0+\ldots+w_3a_3},$$
where $\mu_{d_i}$ denotes the group of $d_i$-th roots of unity with generator
$g_i$; the subgroup
$\mu_d$ of the product is generated by $g = g_0g_1g_2g_3$.
Then, by Corollary 5.5.6 of {\bf [5]}, the dual family ${\cal F}(\Delta^*)$
consists of quotients by $\Pi$ of the family of $\Pi$-invariant
hypersurfaces of degree
$d$ in $\bfp(w)$
$$\sum_{w_0i_0+w_1i_1+w_2i_2+w_3i_3 =
d}a_{i_0i_1i_2i_3}x_0^{i_0}x_1^{i_1}x_2^{i_2}x_3^{i_3} = 0.$$

\medskip
{\bf Example (8.1)}. Let us consider the special case where
$w= (1,1,1,1), d = 4$. The family ${\cal F}(\Delta(w))$ is the family of
quartic
hypersurfaces in $\bfp^3$. Modulo projective transformation this family defines
an open subset of the moduli space
${\bf K}_{<4>}$. The group $\Pi$ is isomorphic to $(\bfz/4)^2$ and its two
generators act by the formula:
$$g_1:(x_0,x_1,x_2,x_3) \to  (x_0,\zeta x_1,x_2,\zeta^3x_3),$$
$$g_2:(x_0,x_1,x_2,x_3) \to  (\zeta x_0,x_1,\zeta^3 x_2,x_3),$$
where $\zeta$ is a primitive 4-th root of unity. The $\Pi$-invariant family of
quartics is
the one-dimensional family
$$V(\lambda) : x_0^4+x_1^4+x_2^4+x_3^4+4\lambda x_0x_1x_2x_3 = 0.$$
The quotient $V(\lambda)/\Pi$ is isomorphic to the surface in $\bfp^4$ given by
the equations:
$$u_0u_1u_2u_3 - u_4^4 = 0,\ \ \ u_0+u_1+u_2+u_3+4\lambda u_4 = 0.\eqno (*)$$
If $\lambda^4 \ne 1$, the surface $V(\lambda)/\Pi$ has six rational double
points of
type $A_3$. Let $V_\lambda$ be the family of K3 surfaces obtained by
simultaneous resolution of
singularities
of the surfaces $V(\lambda)/\Pi, \lambda^4 \ne 1$.

\medskip
\proclaim\ \ \  Theorem (8.2). The family of surfaces $V_\lambda$ is a family
of
$M_2$-polarized surfaces.

{\sl Proof.}
Consider the following four lines on the surface
$$l_i :u_i = u_3 = 0, i= 0,1,2,\ \ l_3:u_0+u_1+u_2 = u_3 = 0.$$
It is easy to check that the six points $P_{ij} = l_i\cap l_j, 0\le i<j\le 3,$
are the singular points of type
$A_3$ of $V(\lambda)/\Pi$. Let $D_{ij}$ be the 6 exceptional divisors coming
from a  minimal
resolution of singularities
$V_\lambda\to V(\lambda)/\Pi$, and let $R_i, i = 0,\ldots,3,$ be the proper
inverse transforms of
the lines. Each divisor $D_{ij}$ consists of three irreducible $(-2)$-curves
with the intersection graph isomorphic to
the Dynkin diagram of type $A_3$. Let $S$ be the sublattice of $\pic(X)$
spanned by the curves
$R_i$ and the irreducible components of the divisors $D_{ij}.$
We shall show that $S\cong M_2$. Consider the divisor
$$D= R_0+R_1+R_2+R_3+R_4+D_{01}+D_{12}+D_{23}+D_{03}.$$
The linear system $|D|$ defines an elliptic fibration on $V_\lambda$ with
reducible
fibre $D$ of type $\tilde A_{15}$ (Kodaira's $I_{16}$). Let $E_{02}$ and
$E_{13}$ be the irreducible
components of the divisors
$D_{02}$ and $D_{13}$ which are disjoint from the divisor $D$. They must be
components of some
reducible fibres of the elliptic fibration.
Since the sublattice of $\pic(V_\lambda)$ generated by irreducible components
of fibres is of rank at most 19, we have only two possibilities.
Either $E_{02},E_{13}$ are components of one fibre of type $\tilde A_3$, or
there exist irreducible curves $E_{01}'$ and $E_{12}'$ such that
$E_{02}+E_{02}'$ and $E_{13}+E_{13}'$ are two fibres of type $\tilde A_1$. In
the first case we find
that Pic$(V_\lambda)$ is of rank $20$. Since the family (*) admits a
degeneration ($\lambda =\infty$) with infinite local monodromy, its image in
the moduli space
${\bf K}_S$ is not a point. Thus for generic $\lambda$, Pic$(V_\lambda)$ is of
rank $19$, and we
have the second possibility. Let $E_0$ be the component of $D_{02}$ which
intersects $R_0$. Then $E_0$ is a section of our fibration, and as such it must
intersect
the fibre $E_{13}+E_{13}'$ at one point. Since it does not intersect $E_{13}$,
it intersects $E_{13}'$
with multiplicity 1. Now we leave to the reader to verify that $V_\lambda$
admits another elliptic fibration which contains
$E_{02}+E_0+E_{13}'+D_{01}$ in its fibre of type $\tilde D_9$ and $D_{23}$ in
its fibre of type $\tilde E_8$.
Arguing as in Example (7.8), we deduce from this that Pic$(X_\lambda) = S \cong
M_2$.

\medskip
 Notice also that the period map
$$p:{\bf A}^1\setminus \{\lambda:\lambda^4 = 1\}\to {\bf K}_{M_2}, \ \  \lambda
\to [X_\lambda]$$
is of degree 4. Indeed, the group ${\bf\mu}_4$ of 4th roots of unity acts on
${\cal F}$ by the formula
$\lambda \to \rho \lambda$ so that $p$ factors through a map
$p':{\bf A}^1\setminus \{1\}\to {\bf A}^1.$
The map $p'$ can also be extended to a map $\bar p':{\bf A}^1\to {\bf A}^1$ by
sending $1$ to the period of a minimal nonsingular model of the
surface $V(1)/\Pi$. This surface represents the unique isomorphism class of
pseudo-ample but not ample
$M_2$-polarized K3 surface. One can show by computing the monodromy at infinity
that
$\bar p'$ is an isomorphism.

\medskip
{\bf Example (8.3)}.
Let $w = (3,1,1,1), d = 6$.
The dual polyhedron $\Delta^*$ can be identified with the convex
hull of the vectors $(1,0,0),(0,1,0),(0,0,1),$
and $(-1,-1,-3).$ The toric hypersurfaces defining the family ${\cal
F}(\Delta^*)$ are given by the
Laurent polynomials
$$aT_1+bT_2+cT_3+dT_1^{-1}T_2^{-1}T_3^{-3}+e = 0.$$
Multiplying both sides by $T_1T_2T_3^3$ and homogenizing, we obtain a
projective model of
$V\in {\cal F}(\Delta^*)$ defined by the equation
$$aT_1^2T_2T_3^3+bT_1T_2^2T_3^3+cT_1T_2T_3^4+dT_0^6 +eT_0T_1T_2T_3^3= 0.$$
This model is not normal. To normalize it, we introduce a new variable $T_4 =
T_0^2/T_3$. Then
a normal projective model can be given by the equations
$$T_1T_2(aT_1+bT_2+cT_3+eT_0)+dT_4^3 = 0,\ \ T_0^2 = T_4T_3.$$
After some obvious linear transformation of the variables, we may assume that
the generic member of
${\cal F}(\Delta(w)^*)$ is isomorphic to the surface $X_\lambda$ in ${\bf P}^4$
given by the equations
$$u_1u_2u_3 - u_4^3 = 0,\ \  (\lambda u_0+u_1+u_2+u_3)u_4+u_0^2 = 0.$$
This is a double cover of the cubic surface $u_1u_2u_3 - u_4^3 = 0$ in ${\bf
P}^3$ branched
along the union of two curves $C_1$ and $C_2$ cut out by the planes $u_4 = 0$
and
$4(u_1+u_2+u_3) -\lambda^2u_4= 0$, respectively. The cubic surface has three
singular points which are cyclic
singularities of type $A_2$. After we resolve them, and then resolve the base
points of the pencil of elliptic curves spanned by
the inverse transforms of the curves $C_1$ and $C_2$, we find a rational
elliptic surface $V_\lambda$
with a singular fibre of type $\tilde A_8$ (originating from the curve $C_1$).
Its double cover branched over this fibre and another
fibre (originating from $C_2$) is birationally isomorphic to $X_\lambda$. After
we resolve its singular points, we obtain
an elliptic K3 surface $\bar X_\lambda$ with a reducible fibre of type $\tilde
A_{17}$. The elliptic
fibration has also three disjoint sections. They come from the three
exceptional curves on
$V_\lambda$
obtained from the resolution of the base points of the elliptic pencil on the
nonsingular model of
the cubic surface. Applying Lemma (7.7), we obtain that $\pic(\bar X_\lambda)$
is a hyperbolic
lattice of rank 19 and discriminant 2. There is only one such lattice, up to
isomorphism.  This is the lattice
$U\perp E_8\perp E_8\perp<-2>.$ Thus
the mirror
family for ${\bf K}_{<2>}$ considered in Example (7.8) can be represented by
the
surfaces from the family ${\cal
F}(\Delta(w)^*)$. Also observe that we have demonstrated the existence of two
different elliptic fibrations on $X$ from the list given in Remark (7.11).

\medskip
{\bf Example (8.4).} In our next example we take
$$w = (1,6,14,26),\ \  d = 42,\  (d_1,d_2,d_3,d_4) =(42,7,3,2).$$
In this case the group $\Pi$ is trivial, and according to Batyrev we should
have
the self-mirrored family. This is true for our mirrors too. The family is
${\bf K}_M$, where
$$M \cong \check M = U\perp E_8.$$
We shall see the latter family again in the next example.

\medskip
{\bf Example (8.5)}.
Consider the affine surface
$$x_1^{d_1}+x_2^{d_2}+x_3^{d_3} = 0,\ \ d_1^{-1}+d_2^{-1}+d_3^{-1}
< 1.$$
According to {\bf [24]} the link space $L$ of the singular point $0$ is
diffeomorphic to the
quotient $G/[\tilde \Gamma,\tilde \Gamma] $, where
$G$ is the universal cover of $PSL(2,\bfr)$ and $[\tilde \Gamma,\tilde \Gamma]$
is the commutator subgroup of
the discrete group $\tilde \Gamma$ of $G$ isomorphic to an extension
$$1\to \bfz \to \tilde \Gamma \to \Gamma(d_1,d_2,d_3)\to 1.$$
Here $\Gamma(d_1,d_2,d_3)$ is the Fuchsian subgroup of $PSL(2,\bfr)$ of
signature $(0;d_1,d_2,d_3).$
Let $K = \tilde \Gamma/[\tilde \Gamma,\tilde \Gamma]$. Its order is
$d_1d_2d_3/d$, where
$d=l.c.m.(d_1,d_2,d_3).$
The quotient
$L/K = G/\tilde \Gamma(d_1,d_2,d_3)$ is the link space of a quasi-homogeneous
triangle
singularity $D_{d_1,d_2,d_3}$(see {\bf [10,29]}). There exist exactly $14$
triples
$d_1,d_2,d_3$
for which
the singularity $D_{d_1,d_2,d_3}$ is isomorphic to the singularity at the
origin of the
affine surface
$P(x,y,z) = 0$, where $P$ is a quasi-homogeneous polynomial of degree $N$ with
weights
$(q_1,q_2,q_3)$ given in the following table:
$$\halign{
# &\hfil\quad#&\hfil\quad#& \hfil\quad#&\hfil\quad#&\hfil\quad#&\hfil\quad#\cr
${\rm
name}$&$(d_1,d_2,d_3)$&$(q_1,q_2,q_3)$&N&$(d_1',d_2',d_3')$&$d_0$&$P(x,y,z)$\cr
$Q_{10}$&$(2,3,9)$&$(6,8,9)$&$24$&$(3,3,4)$&$18$&$z^2x+y^3+x^4$\cr
$Q_{11}$&$(2,4,7)$&$(4,6,7)$&$18$&$(3,3,5)$&$-$&$z^2x+y^3+yx^3$\cr
$Q_{12}$&$(3,3,6)$&$(3,5,6)$&$15$&$(3,3,6)$&6&$z^2x+y^3+x^5$\cr
$Z_{11}$&$(2,3,8)$&$(6,8,15)$&$30$&$(2,4,5)$&24&$y^3x+x^5+z^2$\cr
$Z_{12}$&$(2,4,6)$&$(4,6,11)$&$22$&$(2,4,6)$&12&$y^3x+yx^4+z^2$\cr
$Z_{13}$&$(3,3,5)$&$(3,5,9)$&$18$&$(2,4,7)$&-&$y^3x+x^6+z^2$\cr
$S_{11}$&$(2,5,6)$&$(4,5,6)$&$16$&$(3,4,4)$&-&$z^2x+zy^2+x^4$\cr
$S_{12}$&$(3,4,5)$&$(3,4,5)$&$13$&$(3,4,5)$&-&$y^2z+xz^2+x^3y$\cr
$W_{12}$&$(2,5,5)$&$(4,5,10)$&$20$&$(2,5,5)$&10&$x^5+y^4+z^2$\cr
$W_{13}$&$(3,4,4)$&$(3,4,8)$&$16$&$(2,5,6)$&-&$y^4+yx^4+z^2$\cr
$K_{12}$&$(2,3,7)$&$(6,14,21)$&$42$&$(2,3,7)$&42&$x^7+y^3+z^2$\cr
$K_{13}$&$(2,4,5)$&$(4,10,15)$&$30$&$(2,3,8)$&20&$y^3+yx^5+z^2$\cr
$K_{14}$&$(3,3,4)$&$(3,8,12)$&$24$&$(2,3,9)$&12&$x^8+y^3+z^2$\cr
$U_{12}$&$(4,4,4)$&$(3,4,4)$&$12$&$(4,4,4)$&4&$x^4+y^3+z^3$\cr}$$

For each of the fourteen triples $(d_1,d_2,d_3)$ consider
the family of hypersurfaces of degree $N$ in $\bfp(1,q_1,q_2,q_3)$ given by the
equation
$$Q(w,x,y,z)= P(x,y,z)+\sum a_{ijk} w^{N-iq_1-jq_2-kq_3}x^{i}y^{j}z^k = 0$$
where the monomials $x^{i_1}y^{i_2}z^{i_3}$ form a basis of the Jacobian
algebra
$\bfc[x,y,z]/$(partials of $P$) of the polynomial $P$. There exists a morphism
$Y\to \bfp(1,q_1,q_2,q_3)$
such that the proper inverse transforms of the quasi-smooth hypersurfaces $Q =
0$
form a family ${\cal F}(d_1,d_2,d_3)$  of $M$-polarized K3
surfaces, where
$M$ is the lattice of rank $d_1+d_2+d_3-2$ generated by vectors $e_i$ with
$(e_i,e_i) = -2$ and $(e_i,e_j)\in \{0,1\}, i\ne j,$ determined by the
incidence graph $T_{d_1,d_2,d_3}$ of Dynkin type (for example $T_{2,3,5}$
corresponds to $ E_8$).
Note that the dimension of the family is equal to
${\dim}\ {\bf K}_M$. It is equal to the subscript in the first column minus 2.
The family ${\cal F}(d_1',d_2',d_3')$ corresponds to the mirror moduli space
${\bf K}_{\check M}$.
The involution on the set of fourteen triples
$$(d_1,d_2,d_3) \longleftrightarrow (d_1',d_2',d_3')$$
is the so-called  Arnold's Strange Duality (see {\bf [1,9]}).
If we take the triple $(2,3,7)$ corresponding to the singularity
$K_{12}$ we obtain that ${\cal F}(2,3,7)$ coincides
with the family ${\cal F}(\Delta(w))$ where $w = (1,6,14,21)$. It is self-dual
with respect
to
Batyrev's duality and  mirror duality.

On the other hand let us consider the 12-dimensional family ${\cal F}(3,3,4)$
corresponding to the singularity
$K_{14}$.
It coincides with the family ${\cal F}(\Delta(w))$, where $w = (1,3,8,12).$ The
group $\Pi$ is of order 2.
The Batyrev dual is the 6-dimensional family ${\cal F}(\Delta(w)^*)$ obtained
by dividing $\Pi$-invariant members of
${\cal F}(\Delta(w))$ by $\Pi$. The mirror family of  ${\cal F}(\Delta(w))$
is the 8-dimensional family ${\cal F}(2,3,9)$ corresponding to the singularity
$Q_{10}$.

Suppose $(d_1,d_2,d_3)$ is such that there exists an integer $d_0$ such that
$${1\over d_0}+{1\over d_1}+{1\over d_2}+{1\over d_3} = 1.$$
This happens for 9 triples from the above list.  Then we
can consider the family ${\cal F}(\Delta(w))$ where $w_i = d/d_i, d =
g.c.d.(d_0,d_1,d_2,d_3)$.
The group $\Pi$ is isomorphic to the group $K$ from above. The $\Pi$-quotients
of surfaces from ${\cal F}(\Delta(w))$ are smoothings of the singularity
$D_{d_1,d_2,d_3}$ and hence belong
to the family ${\cal F}(d_1,d_2,d_3)$. This shows that the Batyrev dual family
${\cal F}(\Delta(w)^*)$ is a subfamily of
${\cal F}(d_1,d_2,d_3)$.

For example, ${\cal F}(\Delta(1,3,8,12)^*)$ is a subfamily of ${\cal F}(2,3,8)$
of dimension
9 corresponding to the singularity $Z_{11}$. Also, ${\cal
F}(\Delta(1,4,5,10)^*)$ is a
subfamily of ${\cal F}(2,4,5)$ of
dimension 11 corresponding to the singularity $K_{13}$.

\medskip
The exact relationship between the two mirror constructions seems to be the
following. Let $X$ be a member of the
family ${\cal F}(\Delta)$. Then ${\rm Pic}(X)$ contains the primitive
sublattice
generated by the image of
the restriction homomorphism ${\rm Pic}(\tilde \bfp_\Delta)\to {\rm Pic}(X)$.
Let $M_\Delta$ be the abstract lattice isomorphic to this
lattice. One can show that $M_\Delta\cong {\rm Pic}(X)$ for general member
${\cal F}(\Delta)$ if and
only if $\Delta$ satisfies the following condition: for any 1-dimensional face
$\Gamma$ of $\Delta$,
$$l^*(\Gamma) = l^*(\Gamma^*) = 0,\eqno (*)$$
where $\Gamma^*$ is the dual one-dimensional face of $\Delta^*$, and
$l^*(F)$ denotes the number of integral points in the interior of a face $F$
(see {\bf [2]}).

\medskip
\proclaim \ \ \ Conjecture (8.6).  The lattice $M_\Delta$ always contains a
$1$-admissible
isotropic
vector such that
there exists a primitive embedding
$$M_{\Delta^*} \subset \check M_\Delta.$$
Moreover, the equality takes place if and only if condition (*) is satisfied.

\medskip
This conjecture is confirmed by a result of Batyrev (unpublished) and Kobayashi
{\bf [19]} implying that
$${\rm rank}\   M_\Delta+{\rm rank}\   M_{\Delta^*} \le 20.$$
Also, Kobayashi shows that ${\rm rank}\   M_{\Delta^*} = {\rm rank}\   \check
M_\Delta$ if (*) is
satisfied.
Finally, the conjecture is consistent with the examples from above. If $\Delta
=
\Delta(1,3,8,12)$, we have
$M_\Delta \cong T_{3,3,4}$. Since ${\cal F}(\Delta^*)$ is a subfamily of ${\cal
F}(2,3,8)$ we have
$M_{\Delta^*} \cong T_{2,3,8}$. On the other hand $\check M_\Delta = T_{2,3,9}$
and obviously
$T_{2,3,8}$ embeds naturally in $T_{2,3,9}$. In the second example where
$\Delta = \Delta(1,4,5,10)$ we have
$M_{\Delta^*}\cong T_{2,4,5}$ and $\check M = T_{2,5,5}$.

In some cases when (*) is not satisfied, it is still possible to find a
polyhedron $\Delta'$ satisfying (*)
and
such that ${\cal F}(\Delta')$ is a subfamily of ${\cal F}(\Delta)$. For
example, it is always
possible in the case of the fourteen families from Example (8.5) {\bf [19]}. In
this case, one can verify that
$M_{\Delta'^*} = \check M_{\Delta'}.$

\medskip
{\bf Remark (8.7).} There is a list of weighted projective K3-hypersurfaces
with
Gorenstein singularities
first derived by
Miles Reid (unpublished) and later rediscovered by Yonemura {\bf [42]}). It
consists of 95 families. It contains the family of quartic hypersurfaces and
its mirror family
represented by surfaces of degree 36 in ${\bf P}(7,8,9,12)$ (number 52 in the
list of
Yonemura).

\vglue .2 in
{\bf 9. Other examples.}  Here we consider the examples related to Enriques and
Kummer surfaces.

\medskip
{\bf Example (9.1).} Let $F$ be an Enriques surface, and $p:X\to F$ be its
K3-cover (see {\bf [4,7]}). We have
$$H^2(F,\bfz)/Tors \cong Pic(F)/Tors \cong E_8\perp U.$$
Thus $p^*(\pic(F))$ is a sublattice of $\pic(X)$ isomorphic to $M = E_8(2)\perp
U(2)$ and $X$ acquires a
canonical structure of an $M$-polarized K3 surface. Since $M$ does not contain
vectors $\delta$ with
$(\delta,\delta) = -2$, we can choose $C(M)^+$ to be equal to $V(M)^+$.
Replacing $j: M\to \pic(X)$
by $-j$, if needed, we may assume that
$j(V(M)^+)$ contains the class of an ample divisor $p^*(D)$, where $D$ is an
ample divisor on $F$.
Thus any marked Enriques surface $(F,\phi:H^2(F,\bfz)/Tors \to E_8\perp U)$
defines an ample $M$-polarized K3 surface
$(X,j)$. Conversely, given such $(X,j)$, it defines an involution $\sigma$ on
$H^2(X,\bfz)$ by
setting $\sigma(v) = x$, for any $v\in j(M)$, and $\sigma(v) = -v$, for any
$v\in (j(M))^\perp.$
One can show that any two primitive embeddings of the lattice $M$ on $L$ differ
by an isometry of $L$
(see {\bf [28]}). Thus, we can choose a marking of $\phi:H^2(X,\bfz)\to L$ such
that $j_\phi = j$.
Since the involution $\sigma$ leaves the period $H^{2,0}(X)$ of $X$ unchanged,
by the
Global Torelli Theorem (Corollary (3.2)), there is a unique involution $\tau$
of $X$ such that
$\sigma = \tau^*$. By using the
the Lefschetz fixed-point formula, it is not difficult to show that the set of
fixed points of
$\tau$ is empty (see {\bf [28]}, p.221). Thus $X = F/(\tau)$ is an Enriques
surface together
with a marking defined by descending the isomorphism $j:E_8(2)\perp U(2)\to
\pic(F)$ to the isomorphism
$E_8\perp U\to \pic(X) = \pic(F)^\tau.$
This esablishes a bijective correspondence between
the isomorphism classes of marked Enriques surfaces and isomorphism classes of
ample $M$-polarized K3 surfaces.
In particular, ${\bf K}^a_{E_8(2)\perp U(2)}$ can be viewed as the moduli space
of marked Enriques surfaces.

We may embed $E_8(2)\perp U(2)$ diagonally into
$E_8\perp E_8\perp U\perp U$ to obtain that
$$(E_8(2)\perp U(2))_L^\perp \cong E_8(2)\perp U(2)\perp U.$$
In particular, if we define the mirror lattice by taking $f\in U$, we obtain
$$\check M \cong M.$$
Thus the moduli space ${\bf K}_M$ is its own mirror.
If we take $f$ from $U(2)$ instead, we obtain
$$\check M = U\perp E_8(2).$$
One can show that the moduli space ${\bf K}_{U\perp E_8(2)}$ can be represented
by the family of double covers of the
plane
branched along the union of two cubics.

\bigskip
{\bf Example (9.2).} The mirror family for the family of nonsingular minimal
models of the Kummer
surfaces $X ={\rm Kum(A)}$
associated to principally
polarized abelian surfaces $A$ is the family ${\bf K}_M$, where $M^\perp =
U(2)\perp
U(2)\perp <-4>$. This must be well known but let me give a proof due to J.
Keum.
By Theorem 1.4.14 from {\bf [31]} the embedding $<2> \hookrightarrow
U\perp U\perp U$
is unique. Therefore we may assume that the class $h$ of the polarization of
$A$ is mapped to
$e+f$ where $e,f\in U, (e,e) = (f,f) = 0, (e,f) = 1$. Therefore
$T_A = {\rm Pic}(A)^\perp_{H^2(A,\bfz)} \cong U\perp U\perp <-2>$. On the other
hand  it follows from ({\bf [4]}, Chapter VIII, \S5) that
$T_X = {\rm Pic}(X)^\perp_{H^2(X,\bfz)} \cong T_A(2) \cong U(2)\perp U(2)\perp
<-4>$.

Now if we take $f$ from a copy of $U(2)$ we get
$$\check M = U(2)\perp
<-4> \cong \pmatrix{0&2&2\cr
2&0&2\cr
2&2&0\cr}.$$
The miror family is the moduli space of hypersurfaces of degree
$(2,2,2)$ in $\bfp^1\times \bfp^1\times \bfp^1$. This family is equal to the
family
${\cal F}(\Delta)$ where $\Delta = [-1,1]^3 \subset \bfr^3$. As was shown by
Batyrev the dual family
${\cal F}(\Delta^*)$ is the family ${\bf K}_{M'}$ where $M'^\perp = U\perp
\check M$.
This family of K3 surfaces was studied by C. Peters and J. Stienstra in ${\bf
[33]}$.
If we take
$f\in U$, the dual moduli space of ${\bf K}_{M'}$ is equal to
${\bf K}_{\check M}$. So the Kummer family and the Peters-Stienstra family
share the same mirror family.

\vglue  .3 in

\noindent
{\bf References}

\medskip\noindent
[1] V. Arnold, {\it Critical points of smooth functions}, Proc. I.C.M,
Vancouver,
1974, pp. 18-39.

\smallskip
\noindent
[2] P. Aspinwall, D. Morrison, {\it String theory on K3 surfaces}, 1994,
IASSNS-hep-94/23

\smallskip
\noindent
[3] P. Aspinwall, D. Morrison, {\it Mirror symmetry and the moduli space of K3
surfaces}, (to appear).

\smallskip\noindent
[4] W. Barth, C. Peters, A. Van de Ven, {\it Compact Complex Surfaces},
Ergenbnisse der Mathematik
und ihrer Grenzgebiete, 3. Folge, Band 4, Springer-Verlag, 1984.

\smallskip
\noindent
[5]  V. Batyrev, {\it Dual polyhedra and mirror symmetry for Calabi-Yau
hypersurfaces in toric varieties}, J. Alg. Geometry, 3 (1994), 493-535.

\smallskip
\noindent
[6] C. Borcea, {\it K3 surfaces with involution and mirror pairs of Calabi-Yau
manifolds},
Rider College, preprint

\smallskip\noindent
[7] F. Cossec, I. Dolgachev, {\it Enriques surfaces I}, Birkh\"auser. 1989.

\smallskip\noindent
[8] I. Dolgachev, V. Nikulin, {\it Exceptional singularities of V. I. Arnold
and K3 surfaces}, Proc. USSR Topological Conference in Minsk, 1977.

\smallskip\noindent
[9] I. Dolgachev, {\it Integral quadratic forms:applications to algebraic
geometry},
Sem. Bourbaki, 1982/83, $n^\circ 611$, Asterisque, vol. 105/106, Soc. Math.
France, pp. 251-275.

\smallskip\noindent
[10] I. Dolgachev, {\it On algebraic properties of algebras of automorphic
forms}, in
``Modular Functions in Analysis and Number Theory'', Lect. Notes in
Mathematics and Statistics, vol. 5, Univ. Pittsburgh, 1983, pp. 21-29.

\noindent
[11] {\it Essays on Mirror Symmetry} (ed. S.-T. Yau), Int. Press Co.,
Hong Kong, 1992.

\smallskip\noindent
[12] R. Fricke, {\it Lehrbuch der Algebra}, B. 3, Braunschweig, 1928.

\smallskip\noindent
[13] {\it G\'eometrie des surfaces K3: modules et p\'eriodes}, Ast\'erisque,
vol.
126, Soc. Math. France, 1985.

\smallskip\noindent
[14] A. Giveon, D.-J. Smit, {\it Symmetries of the moduli space of (2,2)
superstring vacua}, Nucl. Phys. B349 (1991), 168-206.

\smallskip\noindent
[15] Ph. Griffiths, {\it Periods of integrals on algebraic manifolds}, I,II,
Amer. J. Math. 90 (1968), 568-626, 805-865.

\smallskip\noindent
[16] Ph. Griffiths, L. Tu, {\it Infinitesimal variation of Hodge structure}, In
`Topics in Transcendental
Algebraic Geometry'', Ann. Math. Studies, vol. 106 , Princeton University
Press, 1984.

\smallskip\noindent
[17] D. James, {\it On Witt's Theorem for unimodular quadratic forms}, Pac. J.
Math.
26:2 (1968), 303-316.

\smallskip\noindent
[18] P.G. Kluit, {\it On the normalizer of $\Gamma_0(N)$}, in ``Modular
Functions of One Variable, V'',
Lect. Notes in Math., vol. 601, Springer, 1977, pp. 239-246.

\smallskip\noindent
[19] M. Kobayashi, {\it Duality of weights, mirror symmetry and Arnold's
strange duality}, 1994, preprint

\smallskip\noindent
[20] J. Lehner, W. Newman, {\it Weierstrass points of $\Gamma_0(n)^*$}, Ann.
Math., 79 (1964),
360-368.

\smallskip\noindent
[21] B. Lian, S.-T. Yau, {\it Arithmetic properties of mirror map and quantum
coupling}, preprint, 1994, hep-th.

\smallskip\noindent
[22] B. Lian, S.-T. Yau, {\it Mirror maps, modular relations and hypergeometric
series II}, preprint, 1994, hep-th.

\smallskip\noindent
[23] E. Martinec, {\it Criticality, catastrophes, and compactifications}, in
``Physics and Mathematics of Strings'', World Scientific, 1990,pp.389-433.

\smallskip\noindent
[24] J. Milnor, {\it On the 3-dimensional Brieskorn manifolds}, in ``Knots,
groups
and 3-manifolds'', Ann. Math. Studies, vol. 84, Princeton Univ. Press 1975, pp.
175-224.

\smallskip\noindent
[25] D. Morrison, {\it On $K3$ surfaces with large Picard number}, Invent.
Math. 75 (1984),
105--121.

\smallskip\noindent
[26] D. Morrison, {\it Mirror symmetry and rational curves on quintic 3-folds:
A
guide for mathematicians}, J. Amer. Math. Soc. 6 (1993), 223-247.

\smallskip\noindent
[27] M. Nagura, K. Sugiyama, {\it Mirror symmetry of the K3 surface}, Intern.
J. Modern Physics A, vol. 10, No 2 (1995), 233-252.

\smallskip\noindent
[28] Y. Namikawa, {\it Periods of Enriques surfaces}, Math. Ann., 270 (1985),
201-222.

\smallskip\noindent
[29] W. Neumann, {\it Abelian covers of quasihomogeneous singularities}, in
``Singularities'', Proc. Symp. Pure Math., vol. 40, Part 2, A.M.S. Providence,
1983, pp. 233-243.

\smallskip\noindent
[30] V. Nikulin, {\it Finite groups of automorphisms of K\"ahler K3 surfaces},
Proc.
Moscow Math. Society, 38(1980), 71-135.

\smallskip\noindent
[31] V. Nikulin, {\it Integral quadratic forms and some of its geometric
applications}, Izv. Akad. Nauk SSSR, Ser. Math. 43 (1979), 103-167.

\smallskip\noindent
[32] V. Nikulin, {\it On rational maps between $K3$ surfaces},in
``Constantin Caratheodory: an international tribute'', Vol. I, II, World Sci.
Publishing. 1991, pp. 964--995.

\smallskip\noindent
[33] C. Peters, J. Stienstra, {\it A pencil of K3 surfaces related to
Ap\'ery's recurrence for $\zeta(3)$
 and Fermi surfaces for potential zero}, in
``Arithmetics of Complex Manifolds'', Lect. Notes in
Math., vol.1399 , Springer-Verlag, 1989.

\smallskip\noindent
[34] H. Pinkham, {\it Singularit\'es exceptionnelles, la dualit\'e \'etrange
d'Arnold et
les surfaces K-3}, C.R. Acad. Sci. Paris, Ser. A-B, 284 (1977), 615-618.

\smallskip\noindent
[35] S.-S. Roan, {\it Mirror symmetry and Arnold's duality}, 1993, preprint
MPI.

\smallskip\noindent
[36] Y. Ruan, G. Tian, {\it A mathematical theory of quantum cohomology}, Math.
Res. Lett. 1,
no. 2 (1994), 269--278.1994.

\smallskip\noindent
[37] F. Scattone, {\it On the compactification of moduli spaces for algebraic
K3
surfaces}, Mem. A.M.S. 70 (1987). No. 374.

\smallskip\noindent
[38] G. Shimura, {\it Introduction to the arithmetic theory of automorphic
functions}, Publ. Math. Soc. Japan, v. 11,
1971.

\smallskip\noindent
[39] H. Sterk, {\it Lattices and K3 surfaces of degree 6}, Lin. Alg. and Appl.,
226-228 (1995), 297-309.

\smallskip\noindent
[40] A. Todorov, {\it  Some ideas from mirror geometry applied to the moduli
space of K3}, preprint.

\smallskip\noindent
[41] C. Voisin, {\it Miroirs et involutions sur les surfaces K3}, in
``Journ\'ees de G\'eom\'etrie Alg\'ebrique d'Orsay'',
vol. 218, Ast\'erisque,
Soc. Math. France, 1993, pp. 273-323.

\smallskip\noindent
[42] T. Yonemura, {\it  Hypersurface simple K3 singularities}, T\^ohoku Math.
J. 42(1990), 351-380.

\bye